\documentclass{article}
\usepackage{amssymb}
\usepackage{amsfonts}
\usepackage{amsmath}

\input{tcilatex}
\begin{document}

\title{On Bernstein processes generated by hierarchies of linear parabolic
systems in $\mathbb{R}^{d}$}
\author{Pierre-A. Vuillermot$^{\ast ,\ast \ast }$ and Jean-C. Zambrini$%
^{\ast \ast }$ \\
UMR-CNRS 7502, Inst. \'{E}lie Cartan de Lorraine, Nancy, France$^{\ast }$\\
Dept. de Matem\'{a}tica, Universidade de Lisboa, Lisboa, Portugal$^{\ast
\ast }$}
\date{}
\maketitle

\begin{abstract}
In this article we investigate the properties of Bernstein processes
generated by infinite hierarchies of forward-backward systems of decoupled
linear deterministic parabolic partial differential equations defined in $%
\mathbb{R}^{d}$, where $d$ is arbitrary. An important feature of those
systems is that the elliptic part of the parabolic operators may be realized
as an unbounded Schr\"{o}dinger operator with compact resolvent in standard $%
L^{2}$-space. The Bernstein processes we are interested in are in general
non-Markovian, may be stationary or non-stationary and are generated by
weighted averages of measures naturally associated with the pure point
spectrum of the operator. We also introduce time-dependent trace-class
operators which possess most of the attributes of density operators in
Quantum Statistical Mechanics, and prove that the statistical averages of
certain bounded self-adjoint observables usually evaluated by means of such
operators coincide with the expectation values of suitable functions of the
underlying processes. In the particular case where the given parabolic
equations involve the Hamiltonian of an isotropic system of quantum harmonic
oscillators, we show that one of the associated processes is identical in
law with the periodic Ornstein-Uhlenbeck process.
\end{abstract}

\section{ Introduction and outline}

Bernstein (or reciprocal) processes constitute a generalization of Markov
processes and have played an increasingly important r\^{o}le in various
areas of mathematics and mathematical physics over the years, particularly
in view of the recent advances in the Monge-Kantorovitch formulation of
Optimal Transport Theory and Stochastic Geometric Mechanics (see, e.g., \cite%
{albeyaza}, \cite{beurling}-\cite{cruzeirozambrini}, \cite{jamison}, \cite%
{lasalle}-\cite{leoroezamb}, \cite{roellythieullen}, \cite{vuizambrini}-\cite%
{zambrinibis} and the many references therein for a history and other works
on the subject, which trace things back to the pioneering works \cite%
{bernstein} and \cite{schroedinger}). As such they may be intrinsically
defined without any reference to partial differential equations, and may
take values in any topological space countable at infinity as was shown in 
\cite{jamison}. However, in this article we restrict ourselves to the
consideration of Bernstein processes generated by certain systems of
parabolic partial differential equations, whose state space is the Euclidean
space $\mathbb{R}^{d}$ endowed with its Borel $\sigma $-algebra $\mathcal{B}%
_{d}$. We begin with the following:

\bigskip

\textbf{Definition 1. }Let $d\in \mathbb{N}^{+}$ and $T\in \left( 0,+\infty
\right) $ be arbitrary. We say the $\mathbb{R}^{d}$-valued process $Z_{\tau
\in \left[ 0,T\right] }$ defined on the complete probability space $\left(
\Omega ,\mathcal{F},\mathbb{P}\right) $ is a Bernstein process if 
\begin{equation}
\mathbb{E}\left( b(Z_{r})\left\vert \mathcal{F}_{s}^{+}\vee \mathcal{F}%
_{t}^{-}\right. \right) =\mathbb{E}\left( b(Z_{r})\left\vert
Z_{s},Z_{t}\right. \right)  \label{condiexpectations}
\end{equation}%
$\mathbb{P}$-almost everywhere for every bounded Borel measurable function $%
b:\mathbb{R}^{d}\mapsto \mathbb{C}$, and for all $r,s,t$ satisfying $r\in
\left( s,t\right) \subset \left[ 0,T\right] $. In (\ref{condiexpectations}),
the $\sigma $-algebras are%
\begin{equation}
\mathcal{F}_{s}^{+}=\sigma \left\{ Z_{\tau }^{-1}\left( F\right) :\tau \leq
s,\text{ }F\in \mathcal{B}_{d}\right\}  \label{pastalgebra}
\end{equation}%
and%
\begin{equation}
\mathcal{F}_{t}^{-}=\sigma \left\{ Z_{\tau }^{-1}\left( F\right) :\tau \geq
t,\text{ }F\in \mathcal{B}_{d}\right\} ,  \label{futurealgebra}
\end{equation}%
where $\mathbb{E}\left( .\left\vert .\right. \right) $ denotes the
conditional expectation on $\left( \Omega ,\mathcal{F},\mathbb{P}\right) $.

\bigskip

This definition obviously extends that of a Markov process in the sense of a
complete independence of the dynamics of $Z_{\tau \in \left[ 0,T\right] }$
within the interval $\left( s,t\right) $ once $Z_{s}$ and $Z_{t}$ are known,
no matter what the behavior of the process is prior to instant $s$ and after
instant $t$. This last property also shows that there are two time
directions coming into play from the outset, since $\mathcal{F}_{s}^{+}$ may
be interpreted as the $\sigma $-algebra gathering all available information
before time $s$ and $\mathcal{F}_{t}^{-}$ as that collecting all available
information after time $t$. It is therefore no surprise that any system of
parabolic partial differential equations susceptible of generating Bernstein
processes should exhibit two time directions, one pointing toward the future
and one toward the past. Accordingly, we introduce below hierarchies of
partial differential equations which we shall define from adjoint parabolic
Cauchy problems of the form%
\begin{align}
\partial _{t}u(\mathsf{x},t)& =\frac{1}{2}\Delta _{\mathsf{x}}u(\mathsf{x}%
,t)-V(\mathsf{x})u(\mathsf{x},t),\text{ \ }(\mathsf{x},t)\in \mathbb{R}%
^{d}\times \left( 0,T\right] ,  \notag \\
u(\mathsf{x},0)& =\varphi _{0}(\mathsf{x}),\text{ \ \ }\mathsf{x}\in \mathbb{%
R}^{d}  \label{cauchyforward}
\end{align}%
and%
\begin{align}
-\partial _{t}v(\mathsf{x},t)& =\frac{1}{2}\Delta _{\mathsf{x}}v(\mathsf{x}%
,t)-V(\mathsf{x})v(\mathsf{x},t),\text{ \ }(\mathsf{x},t)\in \mathbb{R}%
^{d}\times \left[ 0,T\right) ,  \notag \\
v(\mathsf{x,}T)& =\psi _{T}(\mathsf{x}),\text{ \ \ }\mathsf{x}\in \mathbb{R}%
^{d}  \label{cauchybackward}
\end{align}%
where $\Delta _{\mathsf{x}}$ denotes Laplace's operator with respect to the
spatial variable, and where $\varphi _{0}$ and $\psi _{T}$ are real-valued
functions or measures to be specified below. In the sequel we write $%
L^{2}\left( \mathbb{R}^{d}\right) $ and $L^{\infty }\left( \mathbb{R}%
^{d}\right) $ for the usual Lebesgue spaces of all square integrable and
essentially bounded real- or complex-valued functions on $\mathbb{R}^{d}$,
respectively, and $L_{\mathsf{loc}}^{\infty }\mathsf{\left( \mathbb{R}%
^{d}\right) }$ for the local version of $L^{\infty }\left( \mathbb{R}%
^{d}\right) $, without ever distinguishing notationally between the real and
the complex case. It will indeed be clear from the context which case we are
referring to, or else further specifications will be made. Finally we shall
denote by $\left( .,.\right) _{2}$ the inner product in $L^{2}\left( \mathbb{%
R}^{d}\right) $ which we assume to be linear in the first argument, and by $%
\left\Vert .\right\Vert _{2}$ the corresponding norm.

Throughout this article we impose the following hypothesis regarding $V$,
where $\left\vert .\right\vert $ stands for the usual Euclidean norm:

\bigskip

(H) The real-valued function $V$ is bounded from below and satisfies $V\in
L_{\mathsf{loc}}^{\infty }\mathsf{\left( \mathbb{R}^{d}\right) }$ with $V(%
\mathsf{x})\rightarrow +\infty $ as $\left\vert \mathsf{x}\right\vert
\rightarrow +\infty $.

\bigskip

It is well known that Hypothesis (H) allows the self-adjoint realization of
the elliptic operator on the right-hand side of (\ref{cauchyforward})-(\ref%
{cauchybackward}), which is up to a sign the operator associated with the
closure of the quadratic form%
\begin{equation}
\mathsf{Q}\left( f\right) =\frac{1}{2}\sum_{j=1}^{d}\int_{\mathbb{\mathbb{R}}%
^{d}}\mathsf{dx}\left\vert \frac{\partial f(\mathsf{x})}{\partial x_{j}}%
\right\vert ^{2}+\int_{\mathbb{\mathbb{R}}^{d}}\mathsf{dx}V(\mathsf{x}%
)\left\vert f(\mathsf{x})\right\vert ^{2}  \label{quadraticform}
\end{equation}%
first defined for every complex-valued, compactly supported and smooth
function $f$ on $\mathbb{\mathbb{R}}^{d}$ (see, e.g., Section 2 in Chapter
VI of \cite{kato}). Moreover the self-adjoint realization of the operator
associated with (\ref{quadraticform}), henceforth denoted by 
\begin{equation}
H=-\frac{1}{2}\Delta _{\mathsf{x}}+V,  \label{hamiltonian}
\end{equation}%
generates a symmetric semigroup $\exp \left[ -tH\right] $ on $L^{2}\left( 
\mathbb{R}^{d}\right) $ whose integral kernel satisfies%
\begin{equation}
\left\{ 
\begin{array}{c}
g(\mathsf{x},t,\mathsf{y})=g(\mathsf{y},t,\mathsf{x}), \\ 
\\ 
c_{1}t^{-\frac{d}{2}}\exp \left[ -c_{1}^{\ast }\frac{\left\vert \mathsf{x}-%
\mathsf{y}\right\vert ^{2}}{t}\right] \leq g(\mathsf{x},t,\mathsf{y})\leq
c_{2}t^{-\frac{d}{2}}\exp \left[ -c_{2}^{\ast }\frac{\left\vert \mathsf{x}-%
\mathsf{y}\right\vert ^{2}}{t}\right]%
\end{array}%
\right.  \label{greenfunction}
\end{equation}%
for all $\mathsf{x},\mathsf{y\in }\mathbb{\ \mathbb{R}}^{d}\mathbb{\ }$and
every $t\in \left( 0,T\right] $, where $c_{1,2}$ and $c_{1,2}^{\ast }$ are
positive constants (see, e.g., Theorem 1 in \cite{aronsonbis} and its
complete proof in \cite{aronson}). At the same time Hypothesis (H) also
implies that the resolvent of the self-adjoint realization of (\ref%
{hamiltonian}) is compact in $L^{2}\left( \mathbb{R}^{d}\right) $. As a Schr%
\"{o}dinger operator this means that its spectrum $\left( E_{\mathsf{m}%
}\right) _{\mathsf{m}\in \mathbb{N}^{d}}$ is entirely discrete with $E_{%
\mathsf{m}}\rightarrow +\infty $ as $\left\vert \mathsf{m}\right\vert
\rightarrow +\infty $, and that there exists an orthonormal basis $\left( 
\mathsf{f}_{\mathsf{m}}\right) _{\mathsf{m}\in \mathbb{N}^{d}}$ consisting
entirely of its eigenfunctions which we shall assume to be real (see, e.g.,
Section XIII.14 in \cite{reedsimon}, which allows for more general
conditions on $V$). In the sequel we shall refer to the function $g$ in (\ref%
{greenfunction}) as the (parabolic) Green function associated with (\ref%
{cauchyforward})-(\ref{cauchybackward}), also called fundamental solution to
(\ref{cauchyforward}) in references \cite{aronsonbis} and \cite{aronson}.

In the context of this article we also assume that%
\begin{equation}
\mathcal{Z}(t):=\sum_{\mathsf{m}\in \mathbb{N}^{d}}\exp \left[ -tE_{\mathsf{m%
}}\right] <+\infty  \label{convergence}
\end{equation}%
for every $t\in \left( 0,T\right] $, so that the strong convergence of 
\begin{equation}
g(\mathsf{x},t,\mathsf{y})=\sum_{\mathsf{m}\in \mathbb{N}^{d}}\exp \left[
-tE_{\mathsf{m}}\right] \mathsf{f}_{\mathsf{m}}(\mathsf{x})\mathsf{f}_{%
\mathsf{m}}(\mathsf{y})  \label{expansion}
\end{equation}%
holds in $L^{2}\left( \mathbb{R}^{d}\times \mathbb{R}^{d}\right) $. Then the
construction of Bernstein processes associated with (\ref{cauchyforward})-(%
\ref{cauchybackward}) rests on the availability of Green's function (\ref%
{greenfunction}) and on the existence of probability measures on $\mathcal{B}%
_{d}\times \mathcal{B}_{d}$ whose joint densities $\mu $ satisfy the
normalization condition%
\begin{equation}
\int_{\mathbb{R}^{d}\times \mathbb{R}^{d}}\mathsf{dxdy}\mu (\mathsf{x,y})=1.
\label{normalization}
\end{equation}%
Given these facts we organize the remaining part of this article in the
following way: in Section 2 we use the knowledge of $g$ and $\mu $ to state
a general proposition about the existence of a probability space which
supports a Bernstein process $Z_{\tau \in \left[ 0,T\right] }$ characterized
by its finite-dimensional distributions, the joint distribution of $Z_{0}$
and $Z_{T}$ and the probability of finding $Z_{t}$ at any time $t\in \left[
0,T\right] $ in a given region of space. In Section 3 we proceed with the
construction of specific families of probability measures by introducing the
hierarchies of equations we alluded to above. That is, with each level $%
\mathsf{m}$ of the spectrum of (\ref{hamiltonian}) we associate a pair of
adjoint Cauchy problems of the form%
\begin{align}
\partial _{t}u(\mathsf{x},t)& =\frac{1}{2}\Delta _{\mathsf{x}}u(\mathsf{x}%
,t)-V(\mathsf{x})u(\mathsf{x},t),\text{ \ }(\mathsf{x},t)\in \mathbb{R}%
^{d}\times \left( 0,T\right] ,  \notag \\
u(\mathsf{x},0)& =\varphi _{\mathsf{m,0}}(\mathsf{x}),\text{ \ \ }\mathsf{x}%
\in \mathbb{R}^{d}  \label{cauchyforwardbis}
\end{align}%
and%
\begin{align}
-\partial _{t}v(\mathsf{x},t)& =\frac{1}{2}\Delta _{\mathsf{x}}v(\mathsf{x}%
,t)-V(\mathsf{x})v(\mathsf{x},t),\text{ \ }(\mathsf{x},t)\in \mathbb{R}%
^{d}\times \left[ 0,T\right) ,  \notag \\
v(\mathsf{x},T)& =\psi _{\mathsf{m,}T}(\mathsf{x}),\text{ \ \ }\mathsf{x}\in 
\mathbb{R}^{d}.  \label{cauchybackwardbis}
\end{align}%
To wit, we are considering as many pairs of such systems as is necessary to
take into account the whole pure point spectrum of (\ref{hamiltonian}), and
then focus our attention on the sequence of probability measures $\mu _{%
\mathsf{m}}$ given by the joint densities%
\begin{equation}
\mu _{\mathsf{m}}(\mathsf{x,y})=\varphi _{\mathsf{m,0}}(\mathsf{x)}g(\mathsf{%
x},T,\mathsf{y})\psi _{\mathsf{m,}T}(\mathsf{y})  \label{densitybis}
\end{equation}%
where%
\begin{equation}
\left\{ 
\begin{array}{c}
\varphi _{\mathsf{m,0}}(\mathsf{x})=\frac{\delta \left( \mathsf{x}-\mathsf{a}%
_{\mathsf{m}}\right) }{g^{\frac{1}{2}}(\mathsf{a}_{\mathsf{m}},T,\mathsf{b}_{%
\mathsf{m}})}, \\ 
\\ 
\psi _{\mathsf{m,}T}(\mathsf{x})=\frac{\delta \left( \mathsf{x}-\mathsf{b}_{%
\mathsf{m}}\right) }{g^{\frac{1}{2}}(\mathsf{a}_{\mathsf{m}},T,\mathsf{b}_{%
\mathsf{m}})},%
\end{array}%
\right.  \label{inifinacond}
\end{equation}%
thus having%
\begin{equation}
\mu _{\mathsf{m}}(G)=\int_{G}\mathsf{dxdy}\varphi _{\mathsf{m,0}}(\mathsf{x)}%
g(\mathsf{x},T,\mathsf{y})\psi _{\mathsf{m,}T}(\mathsf{y})
\label{positivemeasures}
\end{equation}%
for every $G\in \mathcal{B}_{d}\times \mathcal{B}_{d}$. In the preceding
expressions the points $\mathsf{a}_{\mathsf{m}},\mathsf{b}_{\mathsf{m}}\in 
\mathbb{R}^{d}$ are arbitrarily chosen for every $\mathsf{m\in }\mathbb{N}%
^{d}$ and $\delta $ stands for the Dirac measure so that%
\begin{equation}
\int_{\mathbb{R}^{d}\times \mathbb{R}^{d}}\mathsf{dxdy}\mu _{\mathsf{m}}(%
\mathsf{x,y})=1,  \label{normalizationbis}
\end{equation}%
in agreement with (\ref{normalization}). In this manner and by applying the
general proposition of Section 2 we obtain a sequence of Markovian bridges $%
Z_{\tau \in \left[ 0,T\right] }^{\mathsf{m}}$ whose properties we analyze
thoroughly. With each level of the spectrum we then associate a weight $p_{%
\mathsf{m}}$ and consider probability measures of the form 
\begin{equation}
\bar{\mu}=\sum_{\mathsf{m\in }\mathbb{N}^{d}}p_{\mathsf{m}}\mu _{\mathsf{m}}%
\text{, \ \ }p_{\mathsf{m}}>0\text{, \ \ }\sum_{\mathsf{m\in }\mathbb{N}%
^{d}}p_{\mathsf{m}}=1,  \label{probabilities}
\end{equation}%
that is, statistical mixtures of the measures $\mu _{\mathsf{m}}$. Yet
another application of the proposition of Section 2 then allows us to
generate a non-stationary and non-Markovian process $\bar{Z}_{\tau \in \left[
0,T\right] }$ associated with $\bar{\mu}$. We also introduce a linear,
time-dependent trace-class operator which plays the same r\^{o}le as the
so-called density operator in Quantum Statistical Mechanics (see, e.g, \cite%
{vonneumann}), and prove that the statistical averages of certain bounded
self-adjoint observables evaluated by means of that operator coincide with
the expectations of suitable functions of $\bar{Z}_{\tau \in \left[ 0,T%
\right] }$. In Section 4, keeping the same notation as in Section 3 for the
initial-final data in (\ref{cauchyforwardbis}) and (\ref{cauchybackwardbis}%
), we carry out a similar construction but this time with $\varphi _{\mathsf{%
m,0}}$ and $\exp \left[ -TH\right] \psi _{\mathsf{m,}T}$ forming a complete
biorthonormal system in $L^{2}\left( \mathbb{R}^{d}\right) $, thus
satisfying in particular%
\begin{equation}
\left( \varphi _{\mathsf{m,0}},\exp \left[ -TH\right] \psi _{\mathsf{n,}%
T}\right) _{2}=\delta _{\mathsf{m,n}}  \label{biorthogonality}
\end{equation}%
for all $\mathsf{m,n}\in \mathbb{N}^{d}$ where $\exp \left[ -TH\right] $
stands for the Schr\"{o}dinger semigroup generated by (\ref{hamiltonian})
evaluated at the terminal time $t=T$. The simplest system of this kind is

\begin{equation}
\left\{ 
\begin{array}{c}
\varphi _{\mathsf{m,0}}(\mathsf{x})=\mathsf{f}_{\mathsf{m}}(\mathsf{x}), \\ 
\\ 
\psi _{\mathsf{m,}T}(\mathsf{x})=\exp \left[ TE_{\mathsf{m}}\right] \mathsf{f%
}_{\mathsf{m}}(\mathsf{x})%
\end{array}%
\right.  \label{inifinacondbis}
\end{equation}%
where $E_{\mathsf{m}}$ and $\mathsf{f}_{\mathsf{m}}$ are the eigenvalues and
the eigenfunctions introduced above, respectively, but generally speaking a
pair of initial-final data satisfying (\ref{biorthogonality}) always exists
provided $\exp \left[ -TH\right] \psi _{\mathsf{m,}T}$ is sufficiently close
to\textsf{\ }$\mathsf{f}_{\mathsf{m}}$ for every $\mathsf{m\in }\mathbb{N}%
^{d}$ in some sense. This statement essentially comes from an adaptation of
a result by Paley and Wiener according to Theorem XXXVII of Chapter VII in 
\cite{paleywiener}, but then the corresponding measures (\ref%
{positivemeasures}) are signed since we impose no requirement about the
positivity of $\varphi _{\mathsf{m,0}}$ and $\exp \left[ -TH\right] \psi _{%
\mathsf{m,}T}$. In particular, regarding (\ref{inifinacondbis}) the
eigenfunctions $\mathsf{f}_{\mathsf{m}}$ are typically not positive on $%
\mathbb{R}^{d}$ with the possible exception of $\mathsf{f}_{0}$, so that it
becomes intrinsically impossible to construct a Bernstein process from each $%
\mu _{\mathsf{m}}$ individually in contrast to the method of Section 3.
Nevertheless, the averaging procedure defined by (\ref{probabilities}) still
allows us to generate genuine probability measures on $\mathcal{B}_{d}\times 
\mathcal{B}_{d}$ and hence other non-Markovian processes, which turns out to
be particularly simple to do in the case of (\ref{inifinacondbis}) when%
\begin{equation}
p_{\mathsf{m}}=\mathcal{Z}^{-1}(T)\exp \left[ -TE_{\mathsf{m}}\right]
\label{probagibbs}
\end{equation}%
where $\mathcal{Z}\left( T\right) $ is given by (\ref{convergence}). In
Section 4 we also define a linear, time-dependent trace-class operator from
a pair of suitably chosen Riesz bases and prove again that the corresponding
statistical averages of certain bounded self-adjoint observables coincide
with the expectations of suitable functions of the processes, along with
many other properties. We devote Section 5 to the application of the results
of Sections 3 and 4 to the case where the operator on the right-hand side of
(\ref{cauchyforwardbis})-(\ref{cauchybackwardbis}) is that of an isotropic
system of quantum harmonic oscillators, up to a sign. That is, we consider
hierarchies of the form

\begin{align}
\partial _{t}u(\mathsf{x},t)& =\frac{1}{2}\Delta _{\mathsf{x}}u(\mathsf{x}%
,t)-\frac{\lambda ^{2}}{2}\left\vert \mathsf{x}\right\vert ^{2}u(\mathsf{x}%
,t),\text{ \ }(\mathsf{x},t)\in \mathbb{R}^{d}\times \left( 0,T\right] , 
\notag \\
u(\mathsf{x},0)& =\varphi _{\mathsf{m,0,\lambda }}(\mathsf{x}),\text{ \ \ }%
\mathsf{x}\in \mathbb{R}^{d}  \label{cauchyforwardter}
\end{align}%
and%
\begin{align}
-\partial _{t}v(\mathsf{x},t)& =\frac{1}{2}\Delta _{\mathsf{x}}v(\mathsf{x}%
,t)-\frac{\lambda ^{2}}{2}\left\vert \mathsf{x}\right\vert ^{2}v(\mathsf{x}%
,t),\text{ \ }(\mathsf{x},t)\in \mathbb{R}^{d}\times \left[ 0,T\right) , 
\notag \\
v(\mathsf{x},T)& =\psi _{\mathsf{m,}T,\lambda }(\mathsf{x}),\text{ \ \ }%
\mathsf{x}\in \mathbb{R}^{d}  \label{cauchybackwardter}
\end{align}%
with $\lambda >0$ and suitable choices of $\varphi _{\mathsf{m,0,\lambda }}$
and $\psi _{\mathsf{m,}T,\lambda }$, and prove that the processes
constructed there are intimately tied up with various types of conditioned
Ornstein-Uhlenbeck processes. In particular, we show that one of these is
identical in law with the periodic Ornstein-Uhlenbeck process, which was
also analyzed by means of completely different techniques by various authors
in different contexts (see, e.g., \cite{kwakernaak}, \cite{pedersen} and 
\cite{roellythieullen}). To this end we carry out explicit computations of
the laws and of the covariances based on the fact that in this situation
Green's function identifies with Mehler's $d$-dimensional kernel, namely,%
\begin{eqnarray}
&&g_{\lambda }(\mathsf{x},t,\mathsf{y})  \notag \\
&=&\left( \frac{2\pi \sinh \left( \lambda t\right) }{\lambda }\right) ^{-%
\frac{d}{2}}\exp \left[ -\frac{\lambda \left( \cosh (\lambda t)\left(
\left\vert \mathsf{x}\right\vert ^{2}+\left\vert \mathsf{y}\right\vert
^{2}\right) -2\left( \mathsf{x,y}\right) _{\mathbb{R}^{d}}\right) }{2\sinh
\left( \lambda t\right) }\right]  \label{mehler}
\end{eqnarray}%
for all $\mathsf{x},\mathsf{y\in }\mathbb{R}^{d}$ and every $t\in \left( 0,T%
\right] $, where $\left( .\mathsf{,.}\right) _{\mathbb{R}^{d}}$ stands for
the usual inner product in $\mathbb{R}^{d}$. Finally, we point out that the
periodic Ornstein-Uhlenbeck process we just alluded to has the same law as
one of the processes used in \cite{hoeghkrohn} to discuss properties of
certain quantum systems in equilibrium with a thermal bath, and that it also
identifies with the process "indexed by the circle" and possessing the
"two-sided Markov property on the circle" investigated in \cite{kleinlandau}%
. Our work indeed shows that many of the processes investigated in those
references may be viewed as belonging to a very special class of
non-Markovian and stationary Bernstein processes.

\section{On the existence of Bernstein processes in $\mathbb{R}^{d}$}

Aside from a probability measure $\mu $ on $\mathcal{B}_{d}\times \mathcal{B}%
_{d}$ that satisfies (\ref{normalization}), the typical construction of a
Bernstein process requires a transition function as is the case for Markov
processes. Since there are two time directions provided by (\ref%
{cauchyforward})-(\ref{cauchybackward}) we shall see that the natural choice
is%
\begin{equation}
Q\left( \mathsf{x},t;F,r;\mathsf{y},s\right) :=\dint\limits_{F}\mathsf{dz}%
q\left( \mathsf{x},t;\mathsf{z},r;\mathsf{y},s\right)
\label{transitionfunction}
\end{equation}%
for every $F\in \mathcal{B}_{d}$, where%
\begin{equation}
q\left( \mathsf{x},t;\mathsf{z},r;\mathsf{y},s\right) :=\frac{g(\mathsf{x}%
,t-r,\mathsf{z})g(\mathsf{z},r-s,\mathsf{y})}{g(\mathsf{x},t-s,\mathsf{y})}.
\label{transitiondensity}
\end{equation}%
Both functions are well defined and positive for all $\mathsf{x},\mathsf{y},%
\mathsf{z}\in \mathbb{R}^{d}$ and all $r,s,t$ satisfying $r\in \left(
s,t\right) \subset \left[ 0,T\right] $ by virtue of (\ref{greenfunction}),
and moreover the normalization condition%
\begin{equation*}
Q\left( \mathsf{x},t;\mathbb{R}^{d},r;\mathsf{y},s\right) =1
\end{equation*}%
holds as a consequence of the semigroup composition law for $g$. It is the
knowledge of both $\mu $ and $Q$ that makes it possible to associate a
Bernstein process with (\ref{cauchyforward})-(\ref{cauchybackward}) in the
following way:

\bigskip

\textbf{Proposition 1.} \textit{Let }$\mu $\textit{\ satisfy (\ref%
{normalization}) and let }$Q$ \textit{be given by (\ref{transitionfunction}%
). Then there exists a probability space }$\left( \Omega ,\mathcal{F},%
\mathbb{P}_{\mu }\right) $ \textit{supporting an }$\mathbb{R}^{d}$\textit{%
-valued Bernstein process }$Z_{\tau \in \left[ 0,T\right] }$\textit{\ such
that the following properties are valid:}

\textit{(a) The function }$Q$\textit{\ is the two-sided transition function
of }$Z_{\tau \in \left[ 0,T\right] }$\textit{\ in the sense that} 
\begin{equation}
\mathbb{P}_{\mu }\left( Z_{r}\in F\left\vert Z_{s},Z_{t}\right. \right)
=Q\left( Z_{t},t;F,r;Z_{s},s\right)  \label{transition}
\end{equation}%
\textit{for each }$F\in \mathcal{B}_{d}$ \textit{and all }$r,s,t$\textit{\
satisfying }$r\in \left( s,t\right) \subset \left[ 0,T\right] $. \textit{%
Moreover,}%
\begin{equation}
\mathbb{P}_{\mu }\left( Z_{0}\in F_{0},Z_{T}\in F_{T}\right) =\mu \left(
F_{0}\times F_{T}\right)  \label{jointdistribution}
\end{equation}%
\textit{for all }$F_{0},F_{T}\in \mathcal{B}_{d}$\textit{, that is, }$\mu $%
\textit{\ is the joint probability distribution of }$Z_{0}$\textit{\ and }$%
Z_{T}$\textit{. In particular we have}%
\begin{equation}
\mathbb{P}_{\mu }\left( Z_{0}\in F\right) =\mu \left( F\times \mathbb{R}%
^{d}\right)  \label{probability2}
\end{equation}%
\textit{and}%
\begin{equation}
\mathbb{P}_{\mu }\left( Z_{T}\in F\right) =\mu \left( \mathbb{R}^{d}\times
F\right)  \label{probability3}
\end{equation}%
\textit{for each} $F\in \mathcal{B}_{d}$.

\textit{(b) For every }$n\in \mathbb{N}^{+}$ \textit{the finite-dimensional
distributions of the process are given by}%
\begin{eqnarray}
&&\mathbb{P}_{\mu }\left( Z_{t_{1}}\in F_{1},...,Z_{t_{n}}\in F_{n}\right) 
\notag \\
&=&\int_{\mathbb{R}^{d}\times \mathbb{R}^{d}}\frac{\mathsf{d}\mu \mathsf{%
\left( \mathsf{x,y}\right) }}{g(\mathsf{x},T,\mathsf{y})}\int_{F_{1}}\mathsf{%
dx}_{1}...\int_{F_{n}}\mathsf{dx}_{n}  \notag \\
&&\times \dprod\limits_{k=1}^{n}g\left( \mathsf{x}_{k},t_{k}-t_{k-1},\mathsf{%
x}_{k-1}\right) \times g\left( \mathsf{y},T-t_{n},\mathsf{x}_{n}\right)
\label{distribution}
\end{eqnarray}%
\textit{for all }$F_{1},...,F_{n}\in \mathcal{B}_{d}$\textit{\ and all }$%
t_{0}=0<t_{1}<...<t_{n}<T$\textit{, where }$\mathsf{x}_{0}=\mathsf{x}$%
\textit{. In particular we have}%
\begin{eqnarray}
&&\mathbb{P}_{\mu }\left( Z_{t}\in F\right)  \notag \\
&=&\int_{\mathbb{R}^{d}\times \mathbb{R}^{d}}\frac{\mathsf{d}\mu \mathsf{%
\left( \mathsf{x,y}\right) }}{g(\mathsf{x},T,\mathsf{y})}\int_{F}\mathsf{dz}%
g\left( \mathsf{x},t,\mathsf{z}\right) g\left( \mathsf{z},T-t,\mathsf{y}%
\right)  \label{probability1}
\end{eqnarray}%
\textit{for each }$F\in \mathcal{B}_{d}$\textit{\ and every} $t\in \left(
0,T\right) $.

\textit{(c) }$\mathbb{P}_{\mu }$\textit{\ is the only probability measure
leading to the above properties.}

\bigskip

There already exists a proof of an abstract version of a related statement
in \cite{jamison} as well as a more analytic version of it in \cite%
{vuizambrini}, so that we limit ourselves here to showing how the basic
quantities of interest can be expressed in terms of Green's function (\ref%
{greenfunction}):

\bigskip

\textbf{Proof.} The existence of $\left( \Omega ,\mathcal{F},\mathbb{P}_{\mu
}\right) $ and of $Z_{\tau \in \left[ 0,T\right] }$ is through Kolmogorov's
extension theorem, with the probability $\mathbb{P}_{\mu }$ defined on
cylindical sets by

\begin{eqnarray*}
&&\mathbb{P}_{\mu }\left( Z_{0}\in F_{0},Z_{t_{1}}\in F_{1},...,Z_{t_{n}}\in
F_{n},Z_{T}\in F_{T}\right) \\
&=&\int_{F_{0}\times F_{T}}\mathsf{d}\mu \mathsf{\left( \mathsf{x,y}\right) }%
\int_{F_{1}}\mathsf{dx}_{1}...\int_{F_{n}}\mathsf{dx}_{n}\dprod%
\limits_{k=1}^{n}q\left( y,T;x_{k},t_{k};x_{k-1},t_{k-1}\right)
\end{eqnarray*}%
for all $F_{0},...,F_{T}\in \mathcal{B}_{d}$\ and all $%
t_{0}=0<t_{1}<...<t_{n}<T$, where $\mathsf{x}_{0}=\mathsf{x}$ and $q$ is
given by (\ref{transitiondensity}). Since%
\begin{eqnarray*}
&&\dprod\limits_{k=1}^{n}q\left( y,T;x_{k},t_{k};x_{k-1},t_{k-1}\right) \\
&=&\dprod\limits_{k=1}^{n}\frac{g(\mathsf{y},T-t_{k},\mathsf{x}_{k})g(%
\mathsf{x}_{k},t_{k}-t_{k-1},\mathsf{x}_{k-1})}{g(\mathsf{y},T-t_{k-1},%
\mathsf{x}_{k-1})} \\
&=&\frac{1}{g(\mathsf{x},T,\mathsf{y})}\dprod\limits_{k=1}^{n}g(\mathsf{x}%
_{k},t_{k}-t_{k-1},\mathsf{x}_{k-1})\times g(\mathsf{y},T-t_{n},\mathsf{x}%
_{n})
\end{eqnarray*}%
after $n-1$ cancellations in the products, we therefore obtain%
\begin{eqnarray}
&&\mathbb{P}_{\mu }\left( Z_{0}\in F_{0},Z_{t_{1}}\in F_{1},...,Z_{t_{n}}\in
F_{n},Z_{T}\in F_{T}\right)  \notag \\
&=&\int_{F_{0}\times F_{T}}\frac{\mathsf{d}\mu \mathsf{\left( \mathsf{x,y}%
\right) }}{g(\mathsf{x},T,\mathsf{y})}\int_{F_{1}}\mathsf{dx}%
_{1}...\int_{F_{n}}\mathsf{dx}_{n}  \notag \\
&&\times \dprod\limits_{k=1}^{n}g(\mathsf{x}_{k},t_{k}-t_{k-1},\mathsf{x}%
_{k-1})\times g(\mathsf{y},T-t_{n},\mathsf{x}_{n}),  \label{distributionbis}
\end{eqnarray}%
which is (\ref{distribution}) when $F_{0}=F_{T}=$ $\mathbb{R}^{d}$. We now
prove (\ref{jointdistribution}) by using the symmetry property in (\ref%
{greenfunction}) along with the semigroup composition law for $g$ to get%
\begin{equation}
\int_{\mathbb{R}^{d}}\mathsf{dx}_{1}...\int_{\mathbb{R}^{d}}\mathsf{dx}%
_{n}\dprod\limits_{k=1}^{n}g(\mathsf{x}_{k},t_{k}-t_{k-1},\mathsf{x}%
_{k-1})\times g(\mathsf{y},T-t_{n},\mathsf{x}_{n})=g(\mathsf{x},T,\mathsf{y})
\label{relation}
\end{equation}%
by means of an easy induction argument on $n$. The substitution of (\ref%
{relation}) into (\ref{distributionbis}) with the choice $F_{1}=...=F_{n}=%
\mathbb{R}^{d}$ then leads to the desired relation%
\begin{equation*}
\mathbb{P}_{\mu }\left( Z_{0}\in F_{0},Z_{T}\in F_{T}\right)
=\int_{F_{0}\times F_{T}}\mathsf{d}\mu \mathsf{\left( \mathsf{x,y}\right) ,}
\end{equation*}%
of which (\ref{probability2}) and (\ref{probability3}) are obvious
particular cases. Finally, (\ref{probability1}) is (\ref{distribution}) with 
$n=1$. $\ \ \blacksquare $ \ \ 

\bigskip

\textsc{Remark.} It is plain that the only relevant conditions in the proof
of the proposition are the symmetry and the positivity of Green's function (%
\ref{greenfunction}), aside from the data of a probability measure $\mu $.
Furthermore, as we shall see below Bernstein processes may be stationary and
Markovian but in general they are neither one nor the other, as these
properties are intimately tied up with the structure of $\mu $. More
specifically, according to Theorem 3.1 in \cite{jamison} adapted to the
present situation, a Bernstein process is Markovian if, and only if, there
exist positive measures $\nu _{0}$ and $\nu _{T}$ on $\mathcal{B}_{d}$ such
that $\mu $ be of the form%
\begin{equation}
\mu (G)=\int_{G}\mathsf{d}\left( \nu _{0}\otimes \nu _{T}\right) \left( 
\mathsf{x,y}\right) g(\mathsf{x},T,\mathsf{y})  \label{markovianmeasure}
\end{equation}%
for every $G\in \mathcal{B}_{d}\times \mathcal{B}_{d}$, with%
\begin{equation*}
\mu (\mathbb{R}^{d}\times \mathbb{R}^{d})=1\text{.}
\end{equation*}%
We refer the reader for instance to \cite{vuillermot}, \cite{vuizambrini}
and some of their references for an analysis of the Markovian case in
various situations. Finally, in various different forms Bernstein processes
have also recently appeared in applications of Optimal Transport Theory as
testified for instance in \cite{leonard} and in the monographs \cite%
{galichon} and \cite{villani}, and in the developments of Stochastic
Geometric Mechanics as in \cite{zambrinibis}.

\bigskip

In the next section we carry out the program described in Section 1 starting
with (\ref{cauchyforwardbis}), (\ref{cauchybackwardbis}) and (\ref%
{densitybis}) when the initial-final data are given by (\ref{inifinacond}).

\section{On mixing Bernstein bridges in $\mathbb{R}^{d}$}

Relation (\ref{inifinacond}) implies that measures (\ref{positivemeasures})
are already probability measures so that we may apply Proposition 1 directly
and in this manner associate a Bernstein process $Z_{\tau \in \left[ 0,T%
\right] }^{\mathsf{m}}$ with each level of the spectrum. This leads to the
following result where $u_{\mathsf{m}}$ and $v_{\mathsf{m}}$ denote the
solutions to (\ref{cauchyforwardbis}) and (\ref{cauchybackwardbis}),
respectively, that is,%
\begin{equation}
u_{\mathsf{m}}(\mathsf{x},t)=\int_{\mathbb{R}^{d}}\mathsf{dy}g(\mathsf{x},t,%
\mathsf{y})\varphi _{\mathsf{m,0}}(\mathsf{y})>0  \label{forwardsolution}
\end{equation}%
and%
\begin{equation}
v_{\mathsf{m}}(\mathsf{x},t)=\int_{\mathbb{R}^{d}}\mathsf{dy}g(\mathsf{x}%
,T-t,\mathsf{y})\psi _{\mathsf{m,}T}(\mathsf{y})>0.  \label{backwardsolution}
\end{equation}

\bigskip

\textbf{Theorem 1.} \textit{Assume that Hypothesis }(H)\textit{\ holds. Then
for every }$\mathsf{m}\in \mathbb{N}^{d}$ \textit{there exists a
non-stationary Bernstein process} $Z_{\tau \in \left[ 0,T\right] }^{\mathsf{m%
}}$ \textit{in} $\mathbb{R}^{d}$ \textit{such that the following statements
are valid:}

\textit{(a) The process }$Z_{\tau \in \left[ 0,T\right] }^{\mathsf{m}}$ 
\textit{is a forward Markov process whose finite-dimensional distributions
are}%
\begin{eqnarray}
&&\mathbb{P}_{\mu _{\mathsf{m}}}\left( Z_{t_{1}}^{\mathsf{m}}\in
F_{1},...,Z_{t_{n}}^{\mathsf{m}}\in F_{n}\right)  \notag \\
&=&\int_{\mathbb{R}^{d}}\mathsf{dx}\rho _{\mathsf{m},0}(\mathsf{x}%
)\int_{F_{1}}\mathsf{dx}_{1}...\int_{F_{n}}\mathsf{dx}_{n}\dprod%
\limits_{k=1}^{n}w_{\mathsf{m}}^{\ast }\left( \mathsf{x}_{k-1},t_{k-1};%
\mathsf{x}_{k},t_{k}\right)  \label{distributionter}
\end{eqnarray}%
\textit{for every }$n\in \mathbb{N}^{+}$,\textit{\ all }$F_{1},...,F_{n}\in 
\mathcal{B}_{d}$\textit{\ and all }$0=t_{0}<t_{1}<...<t_{n}<T$\textit{, with 
}$\mathsf{x}_{\mathsf{0}}=\mathsf{x}$. \textit{In the preceding expression
the density of the forward Markov transition function is}%
\begin{equation}
w_{\mathsf{m}}^{\ast }(\mathsf{x},s;\mathsf{y},t)=g(\mathsf{x},t-s,\mathsf{y}%
)\frac{v_{\mathsf{m}}(\mathsf{y},t)}{v_{\mathsf{m}}(\mathsf{x},s)}
\label{markovdensity}
\end{equation}%
\textit{for all} $\mathsf{x,y}\in \mathbb{R}^{d}$ \textit{and all} $s,t\in %
\left[ 0,T\right] $ \textit{with }$t>s$, \textit{while the initial
distribution of the process reads}%
\begin{equation}
\rho _{\mathsf{m},0}(\mathsf{x})=\varphi _{\mathsf{m,0}}(\mathsf{x})v_{%
\mathsf{m}}(\mathsf{x},0).  \label{initialmarginal}
\end{equation}

\textit{(b) The process }$Z_{\tau \in \left[ 0,T\right] }^{\mathsf{m}}$ 
\textit{may also be viewed as a backward Markov process since the
finite-dimensional distributions (\ref{distributionter}) may also be written
as}%
\begin{eqnarray}
&&\mathbb{P}_{\mu _{\mathsf{m}}}\left( Z_{t_{1}}^{\mathsf{m}}\in
F_{1},...,Z_{t_{n}}^{\mathsf{m}}\in F_{n}\right)  \notag \\
&=&\int_{\mathbb{R}^{d}}\mathsf{dx}\rho _{\mathsf{m},T}(\mathsf{x}%
)\int_{F_{1}}\mathsf{dx}_{1}...\int_{F_{n}}\mathsf{dx}_{n}\dprod%
\limits_{k=1}^{n}w_{\mathsf{m}}\left( \mathsf{x}_{k+1},t_{k+1};\mathsf{x}%
_{k},t_{k}\right)  \label{distributionquarto}
\end{eqnarray}%
\textit{for every }$n\in \mathbb{N}^{+}$,\textit{\ all }$F_{1},...,F_{n}\in 
\mathcal{B}_{d}$\textit{\ and all }$0<t_{1}<...<t_{n}<t_{n+1}=T$\textit{,
with }$\mathsf{x}_{n\mathsf{+1}}=\mathsf{x}$. \textit{In the preceding
expression the density of the backward Markov transition function is}%
\begin{equation}
w_{\mathsf{m}}(\mathsf{x},t;\mathsf{y},s)=g(\mathsf{x},t-s,\mathsf{y})\frac{%
u_{\mathsf{m}}(\mathsf{y},s)}{u_{\mathsf{m}}(\mathsf{x},t)}
\label{markovdensitybis}
\end{equation}%
\textit{for all} $\mathsf{x,y}\in \mathbb{R}^{d}$ \textit{and all} $s,t\in %
\left[ 0,T\right] $ \textit{with }$t>s$, \textit{while the final
distribution of the process reads}%
\begin{equation*}
\rho _{\mathsf{m},T}(\mathsf{x})=\psi _{\mathsf{m,}T}(\mathsf{x})u_{\mathsf{m%
}}(\mathsf{x},T).
\end{equation*}

\textit{(c) We have}%
\begin{equation}
\mathbb{P}_{\mu _{\mathsf{m}}}\left( Z_{0}^{\mathsf{m}}=\mathsf{a}_{\mathsf{m%
}}\right) =\mathbb{P}_{\mu _{\mathsf{m}}}\left( Z_{T}^{\mathsf{m}}=\mathsf{b}%
_{\mathsf{m}}\right) =1  \label{probability4}
\end{equation}%
\textit{and}%
\begin{equation}
\mathbb{P}_{\mu _{\mathsf{m}}}\left( Z_{t}^{\mathsf{m}}\in F\right) =\int_{F}%
\mathsf{dx}u_{\mathsf{m}}(\mathsf{x},t)v_{\mathsf{m}}(\mathsf{x},t)
\label{probability5}
\end{equation}%
\textit{for each }$t\in \left( 0,T\right) $ \textit{and every} $F\in 
\mathcal{B}_{d}$.

\textit{(d) Finally,}%
\begin{equation}
\mathbb{E}_{_{\mu _{\mathsf{m}}}}\left( b(Z_{t}^{\mathsf{m}})\right) =\int_{%
\mathbb{R}^{d}}\mathsf{dx}b(\mathsf{x})u_{\mathsf{m}}(\mathsf{x},t)v_{%
\mathsf{m}}(\mathsf{x},t)  \label{expectations}
\end{equation}%
\textit{for each bounded Borel measurable function }$b:\mathbb{R}^{d}\mathbb{%
\mapsto C}$ \textit{and every }$t\in \left( 0,T\right) $.

\bigskip

\textbf{Proof.} From (\ref{markovdensity}) and the semigroup composition law
for $g$ we get 
\begin{equation*}
w_{\mathsf{m}}^{\ast }(\mathsf{x},s;\mathsf{y},t)=\int_{\mathbb{R}^{d}}%
\mathsf{dz}w_{\mathsf{m}}^{\ast }(\mathsf{x},s;\mathsf{z},r)w_{\mathsf{m}%
}^{\ast }(\mathsf{z},r;\mathsf{y},t)
\end{equation*}%
for all $\mathsf{x,y,z}\in \mathbb{R}^{d}$ and every $r\in \left( s,t\right)
\subset \left[ 0,T\right] $, so that the Chapman-Kolmogorov relation%
\begin{equation*}
W_{\mathsf{m}}^{\ast }(\mathsf{x},s;F,t)=\int_{\mathbb{R}^{d}}\mathsf{dy}w_{%
\mathsf{m}}^{\ast }(\mathsf{x},s;\mathsf{y},r)W_{\mathsf{m}}^{\ast }(\mathsf{%
y},r;F,t)
\end{equation*}%
holds for every $F\in \mathcal{B}_{d}$, where%
\begin{equation*}
W_{\mathsf{m}}^{\ast }(\mathsf{x},s;F,t):=\int_{F}\mathsf{dy}w_{\mathsf{m}%
}^{\ast }(\mathsf{x},s;\mathsf{y},t).
\end{equation*}%
Therefore $W_{\mathsf{m}}^{\ast }$ is the transition function of a forward
Markov process with density $w_{\mathsf{m}}^{\ast }$. In order to prove (\ref%
{distributionter}) we start with (\ref{distribution}) into which we
substitute (\ref{densitybis}) to obtain%
\begin{eqnarray}
&&\mathbb{P}_{\mu _{\mathsf{m}}}\left( Z_{t_{1}}^{\mathsf{m}}\in
F_{1},...,Z_{t_{n}}^{\mathsf{m}}\in F_{n}\right)  \label{distributionseptimo}
\\
&=&\int_{\mathbb{R}^{d}}\mathsf{d\mathsf{x}}\varphi _{\mathsf{m,0}}(\mathsf{%
x)}\int_{F_{1}}\mathsf{dx}_{1}...\int_{F_{n}}\mathsf{dx}_{n}\dprod%
\limits_{k=1}^{n}g\left( \mathsf{x}_{k},t_{k}-t_{k-1},\mathsf{x}%
_{k-1}\right) \times v_{\mathsf{m}}(\mathsf{x}_{n},t_{n})  \notag
\end{eqnarray}%
where $\mathsf{x}_{0}=\mathsf{x}$ and $t_{0}=0$. Furthermore, using (\ref%
{markovdensity}) we may rewrite the product in the preceding expression as%
\begin{eqnarray*}
&&\dprod\limits_{k=1}^{n}g\left( \mathsf{x}_{k},t_{k}-t_{k-1},\mathsf{x}%
_{k-1}\right) \times v_{\mathsf{m}}(\mathsf{x}_{n},t_{n}) \\
&=&\dprod\limits_{k=1}^{n}w_{\mathsf{m}}^{\ast }(\mathsf{x}_{k-1},t_{k-1};%
\mathsf{x}_{k},t_{k})\times v_{\mathsf{m}}(\mathsf{x},0)
\end{eqnarray*}%
after $n-1$ cancellations, which eventually leads to Statement (a) by taking
(\ref{initialmarginal}) into account. The proof of Statement (b) is similar
and thereby omitted. Now, from (\ref{densitybis}) and (\ref{probability2})
we have%
\begin{equation}
\mathbb{P}_{\mu _{\mathsf{m}}}\left( Z_{0}^{\mathsf{m}}\in F\right) =\int_{F}%
\mathsf{dx}\varphi _{\mathsf{m,0}}(\mathsf{x})v_{\mathsf{m}}(\mathsf{x}%
,0)=\left\{ 
\begin{array}{c}
0\text{ \ if }\mathsf{a}_{\mathsf{m}}\notin F \\ 
1\text{ \ if }\mathsf{a}_{\mathsf{m}}\in F%
\end{array}%
\right.  \label{probability9}
\end{equation}%
by using the first relation in (\ref{inifinacond}), and similarly from (\ref%
{probability3}) we get%
\begin{equation}
\mathbb{P}_{\mu _{\mathsf{m}}}\left( Z_{T}^{\mathsf{m}}\in F\right) =\int_{F}%
\mathsf{dx}u_{\mathsf{m}}(\mathsf{x},T)\psi _{\mathsf{m,}T}(\mathsf{x}%
)=\left\{ 
\begin{array}{c}
0\text{ \ if }\mathsf{b}_{\mathsf{m}}\notin F \\ 
1\text{ \ if }\mathsf{b}_{\mathsf{m}}\in F%
\end{array}%
\right.  \label{probability10}
\end{equation}%
so that (\ref{probability4}) holds. Moreover, (\ref{probability5}) is an
immediate consequence of (\ref{densitybis}), (\ref{probability1}) and (\ref%
{forwardsolution}), (\ref{backwardsolution}), which proves Statement (c) and
thereby Statement (d). Finally, a glance at (\ref{markovdensity}) shows that
(\ref{distributionter}) lacks translation invariance in time so that $%
Z_{\tau \in \left[ 0,T\right] }^{\mathsf{m}}$ is indeed non-stationary. \ \ $%
\blacksquare $

\bigskip

\textsc{Remarks.} (1) The fact that $Z_{\tau \in \left[ 0,T\right] }^{%
\mathsf{m}}$ is both a forward and a backward Markov process is a
manifestation of its reversibility in the sense of Definition 2 in \cite%
{vuizambrini}, which is also readily apparent in (\ref{probability5}) since
the probability density of finding the process in a given region of space at
a given time is expressed as the product of the forward solution (\ref%
{forwardsolution}) times the backward solution (\ref{backwardsolution}). As
a matter of fact we can also obtain (\ref{probability5}) either from (\ref%
{distributionter}) or from (\ref{distributionquarto}) when $n=1$, and we have%
\begin{equation}
\mathbb{P}_{\mu _{\mathsf{m}}}\left( Z_{t}^{\mathsf{m}}\in \mathbb{R}%
^{d}\right) =\int_{\mathbb{R}^{d}}\mathsf{dx}u_{\mathsf{m}}(\mathsf{x},t)v_{%
\mathsf{m}}(\mathsf{x},t)=1  \label{probability7}
\end{equation}%
for every $t\in \left[ 0,T\right] $ as it should be. Indeed, substituting (%
\ref{inifinacond}) into (\ref{forwardsolution})-(\ref{backwardsolution}) and
the resulting expression into (\ref{probability5}) we get%
\begin{equation*}
\mathbb{P}_{\mu _{\mathsf{m}}}\left( Z_{t}^{\mathsf{m}}\in F\right) =\frac{1%
}{g(\mathsf{a}_{\mathsf{m}},T,\mathsf{b}_{\mathsf{m}})}\int_{F}\mathsf{dx}g(%
\mathsf{a}_{\mathsf{m}},t,\mathsf{x})g(\mathsf{x},T-t,\mathsf{b}_{\mathsf{m}%
})\text{,}
\end{equation*}%
which implies (\ref{probability7}) thanks to the semigroup composition law
for $g$. Finally, we stress the fact that the forward density (\ref%
{markovdensity}) is defined from the backward solution (\ref%
{backwardsolution}), while the backward density (\ref{markovdensitybis}) is
defined from the forward solution (\ref{forwardsolution}), and not the other
way around.

(2) We note that (\ref{distributionseptimo}) may also be written as%
\begin{eqnarray}
&&\mathbb{P}_{\mu _{\mathsf{m}}}\left( Z_{t_{1}}^{\mathsf{m}}\in
F_{1},...,Z_{t_{n}}^{\mathsf{m}}\in F_{n}\right)  \label{distributionoctavo}
\\
&=&\int_{F_{1}}\mathsf{dx}_{1}...\int_{F_{n}}\mathsf{dx}_{n}\dprod%
\limits_{k=2}^{n}g\left( \mathsf{x}_{k},t_{k}-t_{k-1},\mathsf{x}%
_{k-1}\right) \times u_{\mathsf{m}}(\mathsf{x}_{1},t_{1})v_{\mathsf{m}}(%
\mathsf{x}_{n},t_{n})  \notag
\end{eqnarray}%
by carrying out the integral over $\mathsf{x}$ and by taking (\ref%
{forwardsolution}) into account. Relation (\ref{distributionoctavo}) will
play an important r\^{o}le in Section 5 since the integrand determines the
density of the law of the random vector $(Z_{t_{1}}^{\mathsf{m}%
},...,Z_{t_{n}}^{\mathsf{m}})\in \mathbb{R}^{nd}$.

(3) Relation (\ref{probability4}) shows that the process $Z_{\tau \in \left[
0,T\right] }^{\mathsf{m}}$ is pinned down at $\mathsf{a}_{\mathsf{m}}$ when $%
t=0$ and at $\mathsf{b}_{\mathsf{m}}$ when $t=T$. We have therefore obtained
a sequence of Markovian bridges associated with the discrete spectrum of the
operator on the right-hand side of (\ref{cauchyforwardbis})-(\ref%
{cauchybackwardbis}), which we shall call Bernstein bridges in the sequel.
In particular, each process $Z_{\tau \in \left[ 0,T\right] }^{\mathsf{m}}$
reduces to a Markovian loop in $\mathbb{R}^{d}$ when $\mathsf{a}_{\mathsf{m}%
}=\mathsf{b}_{\mathsf{m}}$ for every $\mathsf{m}$.

\bigskip

It turns out that Theorem 1 is the stepping stone toward the construction of
a non-Markovian process we alluded to at the beginning of this article,
which we shall carry out through the averaging procedure briefly sketched in
the introduction. Accordingly, by mixing the Bernstein bridges constructed
above we obtain the following result:

\bigskip

\textbf{Theorem 2. }\textit{Assume that Hypothesis }(H)\textit{\ holds, and
for every }$\mathsf{m}\in \mathbb{N}^{d}$ \textit{let} $Z_{\tau \in \left[
0,T\right] }^{\mathsf{m}}$ \textit{be} \textit{the process of Theorem 1. Let 
}$\bar{Z}_{\tau \in \left[ 0,T\right] }$ \textit{be the Bernstein process in
the sense of Proposition 1 where the probability measure is (\ref%
{probabilities}) with the initial-final conditions given by (\ref%
{inifinacond}). Then the following statements are valid:}

\textit{(a) The process }$\bar{Z}_{\tau \in \left[ 0,T\right] }$ \textit{is
non-stationary, non-Markovian and its finite-dimensional distributions are}%
\begin{eqnarray}
&&\mathbb{P}_{\bar{\mu}}\left( \bar{Z}_{t_{1}}\in F_{1},...,\bar{Z}%
_{t_{n}}\in F_{n}\right)  \notag \\
&=&\sum_{\mathsf{m\in }\mathbb{N}^{d}}p_{\mathsf{m}}\mathbb{P}_{\mu _{%
\mathsf{m}}}\left( Z_{t_{1}}^{\mathsf{m}}\in F_{1},...,Z_{t_{n}}^{\mathsf{m}%
}\in F_{n}\right)  \label{distributionsesto}
\end{eqnarray}%
\textit{for every }$n\in \mathbb{N}^{+}$ \textit{and} \textit{all }$%
F_{1},...,F_{n}\in \mathcal{B}_{d}$, \textit{where} $\mathbb{P}_{\mu _{%
\mathsf{m}}}\left( Z_{t_{1}}^{\mathsf{m}}\in F_{1},...,Z_{t_{n}}^{\mathsf{m}%
}\in F_{n}\right) $ \textit{is given either by (\ref{distributionter}) or (%
\ref{distributionquarto}). }

\textit{(b) We have}%
\begin{equation}
\mathbb{P}_{\bar{\mu}}(\bar{Z}_{t}\in F)=\sum_{\mathsf{m\in }\mathbb{N}%
^{d}}p_{\mathsf{m}}\mathbb{P}_{\mu _{\mathsf{m}}}\left( Z_{t}^{\mathsf{m}%
}\in F\right)  \label{probability8}
\end{equation}%
\textit{for each }$t\in \left[ 0,T\right] $ \textit{and every} $F\in 
\mathcal{B}_{d}$, \textit{where }$\mathbb{P}_{\mu _{\mathsf{m}}}\left(
Z_{t}^{\mathsf{m}}\in F\right) $ \textit{is given by (\ref{probability5}), (%
\ref{probability9}) and (\ref{probability10}). }

\textit{(c) We have}%
\begin{equation}
\mathbb{E}_{_{\bar{\mu}}}(b(\bar{Z}_{t})=\sum_{\mathsf{m\in }\mathbb{N}%
^{d}}p_{\mathsf{m}}\mathbb{E}_{_{\mu _{\mathsf{m}}}}(b(Z_{t}^{\mathsf{m}}))
\label{expectationsbis}
\end{equation}%
\textit{for each bounded Borel measurable function }$b:\mathbb{R}^{d}\mathbb{%
\mapsto C}$ \textit{and every }$t\in \left[ 0,T\right] $, \textit{where} $%
\mathbb{E}_{_{\mu _{\mathsf{m}}}}\left( b(Z_{t}^{\mathsf{m}})\right) $%
\textit{\ is given by (\ref{expectations}).}

\bigskip

\textbf{Proof.} It is sufficient to substitute the joint density%
\begin{equation*}
\bar{\mu}(\mathsf{x,y})=g(\mathsf{x},T,\mathsf{y})\sum_{\mathsf{m}\in 
\mathbb{N}^{d}}p_{\mathsf{m}}\varphi _{\mathsf{m,0}}(\mathsf{x)}\psi _{%
\mathsf{m,}T}(\mathsf{y})
\end{equation*}%
with $\varphi _{\mathsf{m,0}}$ and $\psi _{\mathsf{m,}T}$ given by (\ref%
{inifinacond}) into (\ref{distribution}) and (\ref{probability1}) to obtain (%
\ref{distributionsesto}) and (\ref{probability8}), respectively, from which (%
\ref{expectationsbis}) follows. Owing to the lack of translation invariance
in time of (\ref{distributionquarto}), it is then clear that the process $%
\bar{Z}_{\tau \in \left[ 0,T\right] }$ is also non-stationary. Finally, we
note that $\bar{\mu}$ is not of the form (\ref{markovianmeasure}) so that $%
\bar{Z}_{\tau \in \left[ 0,T\right] }$ is indeed non-Markovian. \ \ $%
\blacksquare $

\bigskip

Having associated an arbitrary weight $p_{\mathsf{m}}$ with each level of
the spectrum of (\ref{hamiltonian}), it is now natural to ask whether there
exists a linear bounded operator $\mathcal{R}\left( t\right) $ acting in $%
L^{2}\left( \mathbb{R}^{d}\right) $ for every $t\in \left( 0,T\right) $
possessing most of the attributes of a so-called density operator in Quantum
Statistical Mechanics. If so, an interesting question is to know whether the
averages of certain bounded self-adjoint observables computed by means of
such a density operator are in one way or another related to some
expectation values of the process $\bar{Z}_{\tau \in \left[ 0,T\right] }$.
We shall see that the answer is affirmative if we define%
\begin{equation}
\mathcal{R}\left( t\right) f:=\sum_{\mathsf{m}\in \mathbb{N}^{d}}p_{\mathsf{m%
}}\left( f,u_{_{\mathsf{m}}}(\mathsf{.},t)\right) _{2}v_{\mathsf{m}}(\mathsf{%
.},t)  \label{timedensityoperator}
\end{equation}%
for each complex-valued $f\in L^{2}\left( \mathbb{R}^{d}\right) $ and every $%
t\in \left( 0,T\right) $, where $u_{\mathsf{m}}(\mathsf{.},t)$ and $v_{%
\mathsf{m}}(\mathsf{.},t)$ are given by%
\begin{equation}
u_{\mathsf{m}}(\mathsf{x},t)=\frac{g(\mathsf{x},t,\mathsf{a}_{\mathsf{m}})}{%
g^{\frac{1}{2}}(\mathsf{a}_{\mathsf{m}},T,\mathsf{b}_{\mathsf{m}})}
\label{forwardsolution bis}
\end{equation}%
and 
\begin{equation}
v_{\mathsf{m}}(\mathsf{x},t)=\frac{g(\mathsf{x},T-t,\mathsf{b}_{\mathsf{m}})%
}{g^{\frac{1}{2}}(\mathsf{a}_{\mathsf{m}},T,\mathsf{b}_{\mathsf{m}})},
\label{backwardsolutionbis}
\end{equation}%
respectively, after substitution of (\ref{inifinacond}) into (\ref%
{forwardsolution}) and (\ref{backwardsolution}). We begin with the following
result in whose proof we write $c$ for all the irrelevant positive constants
depending only on the universal constants $c_{1,2}$ and $c_{1,2}^{\ast }$ in
(\ref{greenfunction}):

\bigskip

\textbf{Theorem 3.} \textit{Let us assume that the sequences of points }$a_{%
\mathsf{m}},b_{\mathsf{m}}$\textit{\ in (\ref{inifinacond}) satisfy}%
\begin{equation}
\sup_{\mathsf{m}\in \mathbb{N}^{d}}\left\vert \mathsf{a}_{\mathsf{m}}-%
\mathsf{b}_{\mathsf{m}}\right\vert <+\infty .  \label{supcondition}
\end{equation}%
\textit{\ Then the following statements hold:}

\textit{(a) Formula (\ref{timedensityoperator}) defines a linear trace-class
operator in }$L^{2}\left( \mathbb{R}^{d}\right) $\textit{\ for every }$t\in
\left( 0,T\right) $\textit{\ and we have}%
\begin{equation*}
\func{Tr}\mathcal{R}\left( t\right) =\sum_{\mathsf{m\in }\mathbb{N}^{d}}p_{%
\mathsf{m}}=1\text{.}
\end{equation*}%
\textit{\ }

\textit{(b) Let us consider the linear bounded self-adjoint multiplication
operator on} $L^{2}\left( \mathbb{R}^{d}\right) $ \textit{given by} $Bf=bf$%
\textit{\ for every complex-valued }$f\in L^{2}\left( \mathbb{R}^{d}\right) $%
,\textit{\ where }$b\in $\textit{\ }$L^{\infty }\left( \mathbb{R}^{d}\right) 
$\textit{\ is real-valued. If }$\bar{Z}_{\tau \in \left[ 0,T\right] }$ 
\textit{denotes the Bernstein process of Theorem 2 then we have}%
\begin{equation}
\func{Tr}\left( \mathcal{R}\left( t\right) B\right) =\mathbb{E}_{_{\bar{\mu}%
}}\left( b(\bar{Z}_{t})\right)  \label{averagesequality}
\end{equation}%
\textit{for every} $t\in \left( 0,T\right) $, \textit{where the right-hand
side of (\ref{averagesequality}) is given by (\ref{expectationsbis}).}

\bigskip

\textbf{Proof.} We first prove that $u_{\mathsf{m}}(\mathsf{.},t),v_{\mathsf{%
m}}(\mathsf{.},t)\in L^{2}\left( \mathbb{R}^{d}\right) $ and that there
exists a constant $c_{\ast }>0$ independent of $\mathsf{m}$ and depending
only on $t,T$ and on the constants in (\ref{greenfunction}) such that%
\begin{eqnarray}
\left\Vert u_{\mathsf{m}}(\mathsf{.},t)\right\Vert _{2} &\leq &c_{\ast
}<+\infty ,  \label{uniformbound1} \\
\left\Vert v_{\mathsf{m}}(\mathsf{.},t)\right\Vert _{2} &\leq &c_{\ast
}<+\infty .  \label{uniformbound2}
\end{eqnarray}%
Indeed, from the right-hand side inequality (\ref{greenfunction}) we have%
\begin{equation*}
\int_{\mathbb{R}^{d}}\mathsf{dx}g^{2}(\mathsf{x},t,\mathsf{a}_{\mathsf{m}%
})\leq ct^{-d}\int_{\mathbb{R}^{d}}\mathsf{dx}\exp \left[ -c\frac{\left\vert 
\mathsf{x}\right\vert ^{2}}{t}\right] =ct^{-\frac{d}{2}}<+\infty
\end{equation*}%
for every $t\in \left( 0,T\right) $ independently of $\mathsf{m}$ by
translation invariance of the integral, and similarly%
\begin{equation*}
\int_{\mathbb{R}^{d}}\mathsf{dx}g^{2}(\mathsf{x},T-t,\mathsf{b}_{\mathsf{m}%
})\leq c\left( T-t\right) ^{-\frac{d}{2}}<+\infty .
\end{equation*}%
On the other hand, from the left-hand side inequality (\ref{greenfunction})
we obtain%
\begin{equation*}
\frac{1}{g(\mathsf{a}_{\mathsf{m}},T,\mathsf{b}_{\mathsf{m}})}\leq cT^{\frac{%
d}{2}}\exp \left[ c\frac{\left\vert \mathsf{a}_{\mathsf{m}}-\mathsf{b}_{%
\mathsf{m}}\right\vert ^{2}}{T}\right]
\end{equation*}%
so that we eventually get%
\begin{equation*}
\left\Vert u_{\mathsf{m}}(\mathsf{.},t)\right\Vert _{2}^{2}\leq c\left( 
\frac{T}{t}\right) ^{\frac{d}{2}}\exp \left[ c\frac{\left\vert \mathsf{a}_{%
\mathsf{m}}-\mathsf{b}_{\mathsf{m}}\right\vert ^{2}}{T}\right] \leq c\left( 
\frac{T}{t}\right) ^{\frac{d}{2}}\exp \left[ \frac{c}{T}\right] :=c_{\ast
}^{2}<+\infty
\end{equation*}%
by virtue of (\ref{forwardsolution bis}) and (\ref{supcondition}). In a
completely similar way we have%
\begin{equation*}
\left\Vert v_{\mathsf{m}}(\mathsf{.},t)\right\Vert _{2}^{2}\leq c\left( 
\frac{T}{T-t}\right) ^{\frac{d}{2}}\exp \left[ c\frac{\left\vert \mathsf{a}_{%
\mathsf{m}}-\mathsf{b}_{\mathsf{m}}\right\vert ^{2}}{T-t}\right] \leq
c_{\ast }^{2}
\end{equation*}%
by changing the value of $c_{\ast }$ if necessary, so that (\ref%
{uniformbound1}) and (\ref{uniformbound2}) hold. Therefore, series (\ref%
{timedensityoperator}) converges strongly in $L^{2}\left( \mathbb{R}%
^{d}\right) $ and defines there a linear bounded operator since%
\begin{equation*}
\sum_{\mathsf{m}\in \mathbb{N}^{d}}p_{\mathsf{m}}\left\vert \left( f,u_{%
\mathsf{m}}(\mathsf{.},t)\right) _{2}\right\vert \left\Vert v_{\mathsf{m}}(%
\mathsf{.},t)\right\Vert _{2}\leq c_{\ast }^{2}\left\Vert f\right\Vert
_{2}<+\infty
\end{equation*}%
for each $f\in L^{2}\left( \mathbb{R}^{d}\right) $ and every $t\in \left(
0,T\right) $. In order to prove that $\mathcal{R}\left( t\right) $ is
trace-class, it is then necessary and sufficient to show that%
\begin{equation}
\sum_{\mathsf{n}\in \mathbb{N}^{d}}\left( \mathcal{R}\left( t\right) \mathsf{%
h}_{\mathsf{n}},\mathsf{h}_{\mathsf{n}}\right) _{2}<+\infty
\label{matrixtrace}
\end{equation}%
for any\textit{\ }orthonormal basis $\left( \mathsf{h}_{\mathsf{n}}\right) _{%
\mathsf{n}\in \mathbb{N}^{d}}$ in $L^{2}\left( \mathbb{R}^{d}\right) $, in
which case (\ref{matrixtrace}) will not depend on the choice of the basis
(see, e.g., Theorem 8.1 in Chapter III of \cite{gohbergkrein}). To this end
let us introduce momentarily the auxiliary function%
\begin{equation*}
A(\mathsf{m},\mathsf{n},t):=p_{\mathsf{m}}\left( \mathsf{h}_{\mathsf{n}},u_{%
\mathsf{m}}(\mathsf{.},t)\right) _{2}\left( v_{\mathsf{m}}(\mathsf{.},t),%
\mathsf{h}_{\mathsf{n}}\right) _{2}
\end{equation*}%
so that%
\begin{equation}
\sum_{\mathsf{m}\in \mathbb{N}^{d}}A(\mathsf{m},\mathsf{n},t)=\left( 
\mathcal{R}\left( t\right) \mathsf{h}_{\mathsf{n}},\mathsf{h}_{\mathsf{n}%
}\right) _{2}  \label{identity1}
\end{equation}%
for any fixed $\mathsf{n}$. Moreover, for any fixed $\mathsf{m}$ we have%
\begin{equation}
\sum_{\mathsf{n}\in \mathbb{N}^{d}}A(\mathsf{m},\mathsf{n},t)=p_{\mathsf{m}%
}\left( u_{\mathsf{m}}(\mathsf{.},t),v_{\mathsf{m}}(\mathsf{.},t)\right)
_{2}=p_{\mathsf{m}}  \label{identity2}
\end{equation}%
by virtue of (\ref{probability7}). In addition, the preceding series
converges absolutely as a consequence of the Cauchy-Schwarz inequality and
estimates (\ref{uniformbound1}), (\ref{uniformbound2}) since for any
positive integers $N_{1},...,N_{d}$ we have successively%
\begin{eqnarray}
&&\sum_{\mathsf{n:0\leq n}_{j}\leq N_{j}}\left\vert A(\mathsf{m},\mathsf{n}%
,t)\right\vert  \notag \\
&\leq &p_{\mathsf{m}}\left( \sum_{\mathsf{n}\in \mathbb{N}^{d}}\left\vert
\left( u_{\mathsf{m}}(\mathsf{.},t),\mathsf{h}_{\mathsf{n}}\right)
_{2}\right\vert ^{2}\right) ^{\frac{1}{2}}\left( \sum_{\mathsf{n}\in \mathbb{%
N}^{d}}\left\vert \left( v_{\mathsf{m}}(\mathsf{.},t),\mathsf{h}_{\mathsf{n}%
}\right) _{2}\right\vert ^{2}\right) ^{\frac{1}{2}}  \notag \\
&=&p_{\mathsf{m}}\left\Vert u_{\mathsf{m}}(\mathsf{.},t)\right\Vert
_{2}\left\Vert v_{\mathsf{m}}(\mathsf{.},t)\right\Vert _{2}\leq c_{\ast
}^{2}p_{\mathsf{m}}  \label{estimate7}
\end{eqnarray}%
for any fixed $\mathsf{m}$ so that 
\begin{equation*}
\sum_{\mathsf{n}\in \mathbb{N}^{d}}\left\vert A(\mathsf{m},\mathsf{n}%
,t)\right\vert <+\infty
\end{equation*}%
since the partial sums of this series remain bounded. Finally, (\ref%
{estimate7}) still implies%
\begin{equation*}
\sum_{\mathsf{m}\in \mathbb{N}^{d}}\sum_{\mathsf{n}\in \mathbb{N}%
^{d}}\left\vert A(\mathsf{m},\mathsf{n},t)\right\vert \leq c_{\ast
}^{2}\sum_{\mathsf{m}\in \mathbb{N}^{d}}p_{\mathsf{m}}=c_{\ast }^{2}<+\infty
.
\end{equation*}%
Therefore the corresponding iterated series are equal (see, e.g., Theorem
8.43 in Chapter 8 of \cite{apostol}), that is,%
\begin{equation*}
\sum_{\mathsf{n}\in \mathbb{N}^{d}}\sum_{\mathsf{m}\in \mathbb{N}^{d}}A(%
\mathsf{m},\mathsf{n},t)=\sum_{\mathsf{m}\in \mathbb{N}^{d}}\sum_{\mathsf{n}%
\in \mathbb{N}^{d}}A(\mathsf{m},\mathsf{n},t)
\end{equation*}%
or, equivalently,%
\begin{equation*}
\func{Tr}\mathcal{R}\left( t\right) :=\sum_{\mathsf{n}\in \mathbb{N}%
^{d}}\left( \mathcal{R}\left( t\right) \mathsf{h}_{\mathsf{n}},\mathsf{h}_{%
\mathsf{n}}\right) _{2}=\sum_{\mathsf{m}\in \mathbb{N}^{d}}p_{\mathsf{m}}=1
\end{equation*}%
according to (\ref{identity1}) and (\ref{identity2}), which is (a). As for
the proof of (b), arguing as above for the computation of the trace we have%
\begin{eqnarray*}
&&\func{Tr}\left( \mathcal{R}\left( t\right) B\right) \\
&=&\sum_{\mathsf{n}\in \mathbb{N}^{d}}\sum_{\mathsf{m}\in \mathbb{N}^{d}}p_{%
\mathsf{m}}\left( \mathsf{h}_{\mathsf{n}},bu_{_{\mathsf{m}}}(\mathsf{.}%
,t)\right) _{2}\left( v_{\mathsf{m}}(\mathsf{.},t),\mathsf{h}_{\mathsf{n}%
}\right) _{2} \\
&=&\sum_{\mathsf{m}\in \mathbb{N}^{d}}p_{\mathsf{m}}\left( bu_{_{\mathsf{m}%
}}(\mathsf{.},t),v_{\mathsf{m}}(\mathsf{.},t)\right) _{2}=\mathbb{E}_{_{\bar{%
\mu}}}\left( b(\bar{Z}_{t})\right)
\end{eqnarray*}%
where the last equality follows from (\ref{expectations}) and (\ref%
{expectationsbis}) (note that $u_{_{\mathsf{m}}}(\mathsf{.},t)$ and $v_{%
\mathsf{m}}(\mathsf{.},t)$ are also real-valued). \ \ $\blacksquare $

\bigskip

\textsc{Remarks.} (1) The preceding considerations show that $\mathcal{R}(t)$
is not self-adjoint in general for it is easily seen that its adjoint is
obtained by swapping the r\^{o}le of (\ref{forwardsolution bis}) and (\ref%
{backwardsolutionbis}), that is,%
\begin{equation*}
\mathcal{R}^{\ast }\left( t\right) f=\sum_{\mathsf{m}\in \mathbb{N}^{d}}p_{%
\mathsf{m}}\left( f,v_{_{\mathsf{m}}}(\mathsf{.},t)\right) _{2}u_{\mathsf{m}%
}(\mathsf{.},t).
\end{equation*}%
Aside from that and in addition to the conclusion of Theorem 3, (\ref%
{timedensityoperator}) possesses most of the properties of a density
operator. For instance, every operator $\mathcal{P}^{\mathsf{m}%
}(t):L^{2}\left( \mathbb{R}^{d}\right) \mapsto L^{2}\left( \mathbb{R}%
^{d}\right) $ defined by%
\begin{equation*}
\mathcal{P}^{\mathsf{m}}(t)f:=\left( f,u_{_{\mathsf{m}}}(\mathsf{.}%
,t)\right) _{2}v_{\mathsf{m}}(\mathsf{.},t)
\end{equation*}%
satisfies $\left( \mathcal{P}^{\mathsf{m}}(t)\right) ^{2}=\mathcal{P}^{%
\mathsf{m}}(t)$ as a consequence of (\ref{probability7}) and thus represents
an oblique projection rather than an orthogonal projection, but (\ref%
{timedensityoperator}) is still a statistical mixture of the $\mathcal{P}^{%
\mathsf{m}}(t)$ obtained by sweeping over the whole spectrum of (\ref%
{hamiltonian}). Moreover, we remark that (\ref{timedensityoperator})
involves both the forward and the backward solutions to (\ref%
{cauchyforwardbis}) and (\ref{cauchybackwardbis}), again in agreement with
the fact that there are two time directions in the theory from the outset.

(2) It is tempting to believe that for \textit{every} linear bounded
selfadjoint operator there exists a real-valued $b\in $\textit{\ }$L^{\infty
}\left( \mathbb{R}^{d}\right) $ such that (\ref{averagesequality}) holds,
since such an operator is unitarily equivalent to a multiplication operator
by the spectral theorem. We defer the general analysis of this question to a
separate publication.

\bigskip

In the next section we carry out the program described in the introduction
when the initial-final data satisfy suitable biorthogonality properties, and
where we keep the same notation $\varphi _{\mathsf{m,0}}$ and $\psi _{%
\mathsf{m,}T}$ for them.

\section{On generating Bernstein processes in $\mathbb{R}^{d}$ by mixing
signed measures}

What we first need lies in the following adaptation of a result by Paley and
Wiener (see the abstract form given in Section 86 of Chapter V in \cite%
{riesznagy} of Theorem XXXVII of Chapter VII in \cite{paleywiener}). We omit
the proof as it is essentially available therein modulo trivial changes and
up to the observation that the equality%
\begin{equation*}
\left( \exp \left[ -tH\right] \varphi _{\mathsf{m,0}},\exp \left[ -\left(
T-t\right) H\right] \psi _{\mathsf{n,}T}\right) _{2}=\left( \varphi _{%
\mathsf{m,0}},\exp \left[ -TH\right] \psi _{\mathsf{n,}T}\right) _{2}
\end{equation*}%
holds for every $t\in \left[ 0,T\right] $ as a consequence of the symmetry
of the semigroup $\exp \left[ -tH\right] $. In the following statement all
functions are supposed to be real-valued with $\left( \mathsf{f}_{\mathsf{m}%
}\right) _{\mathsf{m\in }\mathbb{N}^{d}\text{ \ }}$the orthonormal basis of
Section 1:

\bigskip

\textbf{Proposition 2.} \textit{Let }$\left( \psi _{\mathsf{m,}T}\right) _{%
\mathsf{m\in }\mathbb{N}^{d}}$\textit{\ be an arbitrary sequence in }$%
L^{2}\left( \mathbb{R}^{d}\right) $\textit{\ and let us assume that there
exists }$\theta \in \left[ 0,1\right) $ \textit{such that the estimate}%
\begin{equation}
\left\Vert \dsum\limits_{\mathsf{m}\in \mathsf{I}}\gamma _{_{\mathsf{m}%
}}\left( \mathsf{f}_{\mathsf{m}}-\exp \left[ -TH\right] \psi _{\mathsf{m,}%
T}\right) \right\Vert _{2}\leq \theta \left( \sum_{\mathsf{m}\in \mathsf{I}%
}\left\vert \gamma _{\mathsf{m}}\right\vert ^{2}\right) ^{\frac{1}{2}}
\label{estimatebis}
\end{equation}%
\textit{holds for every sequence }$\left( \gamma _{\mathsf{m}}\right) _{%
\mathsf{m\in }\mathbb{N}^{d}}$\textit{\ of real numbers, where the sums in (%
\ref{estimatebis}) run over the same finite set }$I\subset \mathbb{N}^{d}$ 
\textit{which} \textit{may be chosen arbitrarily. Then there exists a
sequence }$\left( \varphi _{\mathsf{m,}0}\right) _{\mathsf{m\in }\mathbb{N}%
^{d}}\subset L^{2}\left( \mathbb{R}^{d}\right) $\textit{\ such that the
following statements are valid:}

\textit{(a) We have}%
\begin{equation}
\left( \exp \left[ -tH\right] \varphi _{\mathsf{m,0}},\exp \left[ -\left(
T-t\right) H\right] \psi _{\mathsf{n,}T}\right) _{2}=\delta _{\mathsf{m,n}}
\label{biorthogonalitybis}
\end{equation}%
\textit{for every} $t\in \left[ 0,T\right] $ \textit{and the strongly
convergent expansions}%
\begin{equation}
\left\{ 
\begin{array}{c}
f=\dsum\limits_{\mathsf{m\in }\mathbb{N}^{d}}\left( f,\varphi _{\mathsf{m,}%
0}\right) _{2}\exp \left[ -TH\right] \psi _{\mathsf{m,}T}, \\ 
\\ 
f=\dsum\limits_{\mathsf{m\in }\mathbb{N}^{d}}\left( f,\exp \left[ -TH\right]
\psi _{\mathsf{m,}T}\right) _{2}\varphi _{\mathsf{m,}0}%
\end{array}%
\right.  \label{expansionbis}
\end{equation}%
\textit{hold} \textit{for every} $f\in L^{2}\left( \mathbb{R}^{d}\right) $.

\textit{(b) The coefficients in (\ref{expansionbis}) satisfy the estimates}%
\begin{eqnarray}
\left( 1+\theta \right) ^{-1}\left\Vert f\right\Vert _{2} &\leq &\left(
\dsum\limits_{\mathsf{m\in }\mathbb{N}^{d}}\left\vert \left( f,\varphi _{%
\mathsf{m,}0}\right) _{2}\right\vert ^{2}\right) ^{\frac{1}{2}}\leq \left(
1-\theta \right) ^{-1}\left\Vert f\right\Vert _{2},  \label{estimate1} \\
\left( 1-\theta \right) \left\Vert f\right\Vert _{2} &\leq &\left(
\dsum\limits_{\mathsf{m\in }\mathbb{N}^{d}}\left\vert \left( f,\exp \left[
-TH\right] \psi _{\mathsf{m,}T}\right) _{2}\right\vert ^{2}\right) ^{\frac{1%
}{2}}\leq \left( 1+\theta \right) \left\Vert f\right\Vert _{2}.
\label{estimate2}
\end{eqnarray}

\bigskip

Thus the sequences $\left( \exp \left[ -TH\right] \psi _{\mathsf{m,}%
T}\right) _{\mathsf{m\in }\mathbb{N}^{d}}$ and $\left( \varphi _{\mathsf{m,}%
0}\right) _{\mathsf{m\in }\mathbb{N}^{d}}$ constitute Riesz bases of $%
L^{2}\left( \mathbb{R}^{d}\right) $ in the terminology of \cite{gohbergkrein}
and it is plain that (\ref{inifinacondbis}) corresponds to $\theta =0$ in
Proposition 2, in which case (\ref{expansionbis}) reduces to the usual
Fourier expansion of $f$ along the orthonormal basis $\left( \mathsf{f}_{%
\mathsf{m}}\right) _{\mathsf{m\in }\mathbb{N}^{d}}$ and (\ref{estimate1}), (%
\ref{estimate2}) to Parseval's equality. The reason why we have to consider $%
\exp \left[ -TH\right] \psi _{\mathsf{m,}T}$ rather than just $\psi _{%
\mathsf{m,}T}$ lies in Relation (\ref{normalizationter}) of the following
result:

\bigskip

\textbf{Lemma 1.} \textit{Let }$\varphi _{\mathsf{m,0}}$ \textit{and} $\exp %
\left[ -TH\right] \psi _{\mathsf{m,}T}$ \textit{be as in Proposition 2 and
let us again define the density }$\mu _{\mathsf{m}}$ \textit{by} 
\begin{equation*}
\mu _{\mathsf{m}}(\mathsf{x,y})=\varphi _{\mathsf{m,0}}(\mathsf{x)}g(\mathsf{%
x},T,\mathsf{y})\psi _{\mathsf{m,}T}(\mathsf{y})\text{.}
\end{equation*}%
\textit{Then the induced measures }$\mu _{\mathsf{m}}$\textit{\ on }$%
\mathcal{B}_{d}\times \mathcal{B}_{d}$ \textit{are} \textit{signed and we
have}%
\begin{equation}
\mu _{\mathsf{m}}\left( \mathbb{R}^{d}\times \mathbb{R}^{d}\right) =1
\label{normalizationter}
\end{equation}%
\textit{for every} $\mathsf{m\in }\mathbb{N}^{d}$\textit{.}

\bigskip

\textbf{Proof.} The measures are signed since there is no requirement about
the positivity of $\varphi _{\mathsf{m,0}}$ and $\exp \left[ -TH\right] \psi
_{\mathsf{m,}T}$. In particular, regarding (\ref{inifinacondbis}) the
eigenfunctions $\mathsf{f}_{\mathsf{m}}$ are typically not positive on $%
\mathbb{R}^{d}$ with the possible exception of $\mathsf{f}_{0}$. Moreover,
expanding $\psi _{\mathsf{m,}T}$ along the orthonormal basis $\left( \mathsf{%
f}_{\mathsf{m}}\right) _{\mathsf{m\in }\mathbb{N}^{d}}$ we get%
\begin{equation*}
\sum_{\mathsf{k}\in \mathbb{N}^{d}}\exp \left[ -TE_{\mathsf{k}})\right]
\left( \varphi _{\mathsf{m,0}},\mathsf{f}_{\mathsf{k}}\right) _{2}\left(
\psi _{\mathsf{n,}T},\mathsf{f}_{\mathsf{k}}\right) _{2}=\delta _{\mathsf{m,n%
}}
\end{equation*}%
from (\ref{biorthogonality}), and therefore%
\begin{equation*}
\mu _{\mathsf{m}}\left( \mathbb{R}^{d}\times \mathbb{R}^{d}\right) =\sum_{%
\mathsf{k}\in \mathbb{N}^{d}}\exp \left[ -TE_{\mathsf{k}})\right] \left(
\varphi _{\mathsf{m,0}},\mathsf{f}_{\mathsf{k}}\right) _{2}\left( \psi _{%
\mathsf{m,}T},\mathsf{f}_{\mathsf{k}}\right) _{2}=1
\end{equation*}%
by substituting (\ref{expansion}) into (\ref{positivemeasures}). \ \ $%
\blacksquare $

\bigskip

The fact that (\ref{probabilities}) may define a probability measure in the
case under consideration is then ensured by the following result:

\bigskip

\textbf{Lemma 2. }\textit{Let the initial-final conditions form a complete
biorthonormal system} \textit{in the sense of Proposition 2, and let }$\bar{%
\mu}$ \textit{be the measure determined by}%
\begin{equation}
\bar{\mu}(G)=\sum_{\mathsf{m\in }\mathbb{N}^{d}}p_{\mathsf{m}}\mu _{\mathsf{m%
}}(G)  \label{probabilitiesbis}
\end{equation}%
\textit{for every }$G\in \mathcal{B}_{d}\times \mathcal{B}_{d}$\textit{. If }%
\begin{equation}
\left( \mathsf{x,y}\right) \mapsto \sum_{\mathsf{m\in }\mathbb{N}^{d}}p_{%
\mathsf{m}}\varphi _{\mathsf{m,0}}(\mathsf{x)}\psi _{\mathsf{m,}T}(\mathsf{y}%
)  \label{positiveschwartz}
\end{equation}%
\textit{is a positive measure} \textit{on} $\mathbb{R}^{d}\times \mathbb{R}%
^{d}$ \textit{then }$\bar{\mu}$\textit{\ is a probability measure on }$%
\mathcal{B}_{d}\mathcal{\times B}_{d}$.

\bigskip

\textbf{Proof.} We have $\bar{\mu}(\mathbb{R}^{d}\times \mathbb{R}^{d})=1$
because of Lemma 1 and the fact that $\sum_{\mathsf{m\in }\mathbb{N}^{d}}p_{%
\mathsf{m}}=1$. On the other hand, the joint density associated with (\ref%
{probabilitiesbis}) reads%
\begin{equation}
\bar{\mu}\left( \mathsf{x,y}\right) =g(\mathsf{x},T,\mathsf{y})\sum_{\mathsf{%
m\in }\mathbb{N}^{d}}p_{\mathsf{m}}\varphi _{\mathsf{m,0}}(\mathsf{x)}\psi _{%
\mathsf{m,}T}(\mathsf{y})  \label{density}
\end{equation}%
where $g(\mathsf{x},T,\mathsf{y})>0$ according to (\ref{greenfunction}). \ \ 
$\blacksquare $

\bigskip

\textsc{Remark.} It may seem abrupt to assume off-hand that (\ref%
{positiveschwartz}) is positive as a measure. However, an important example
illustrating this situation comes about when the initial-final data are
given by (\ref{inifinacondbis}) and the weights associated with the spectrum
by (\ref{probagibbs}). Indeed, in this case we have%
\begin{equation*}
\mu _{\mathsf{m}}\left( \mathsf{x,y}\right) =\exp \left[ TE_{\mathsf{m}}%
\right] g(\mathsf{x},T,\mathsf{y})\mathsf{f}_{\mathsf{m}}(\mathsf{x})\mathsf{%
f}_{\mathsf{m}}(\mathsf{y})
\end{equation*}%
and therefore%
\begin{eqnarray}
\bar{\mu}\left( \mathsf{x,y}\right) &=&\mathcal{Z}^{-1}(T)g(\mathsf{x},T,%
\mathsf{y})\sum_{\mathsf{m\in }\mathbb{N}^{d}}\mathsf{f}_{\mathsf{m}}(%
\mathsf{x})\mathsf{f}_{\mathsf{m}}(\mathsf{y})  \notag \\
&=&\mathcal{Z}^{-1}(T)g(\mathsf{x},T,\mathsf{y})\delta (\mathsf{x}-\mathsf{y}%
)  \label{positivedensity}
\end{eqnarray}%
where the last equality is a consequence of the completeness of the
orthogonal system $\left( \mathsf{f}_{\mathsf{m}}\right) _{\mathsf{m}\in 
\mathbb{N}^{d}}$. It is (\ref{positivedensity}) that will allow us to relate
the above considerations to the periodic Ornstein-Uhlenbeck process in the
next section.

\bigskip

Since the solutions $u_{\mathsf{m}}$ and $v_{\mathsf{m}}$ to (\ref%
{cauchyforwardbis}) and (\ref{cauchybackwardbis}) may now be written in
terms of the Schr\"{o}dinger semigroup defined in Section 1, namely,%
\begin{equation}
u_{\mathsf{m}}(\mathsf{.},t)=\exp \left[ -tH\right] \varphi _{\mathsf{m,0}}
\label{semigroup1}
\end{equation}%
and%
\begin{equation}
v_{\mathsf{m}}(\mathsf{.},t)=\exp \left[ -(T-t)H\right] \psi _{\mathsf{m,}T},
\label{semigroup2}
\end{equation}%
respectively, then by mixing the measures $\mu _{\mathsf{m}}$ as in Lemma 2
we obtain:

\bigskip

\textbf{Theorem 4. }\textit{Assume that Hypothesis }(H)\ \textit{holds, and
let }$\bar{Z}_{\tau \in \left[ 0,T\right] }$ \textit{be the Bernstein
process in the sense of Proposition 1 with }$\bar{\mu}$\textit{\ given by (%
\ref{density}), the particular case (\ref{positivedensity}) being excluded.
Then the following statements are valid:}

\textit{(a) The process }$\bar{Z}_{\tau \in \left[ 0,T\right] }$ \textit{is
non-stationary, non-Markovian and for every }$n\in \mathbb{N}^{+}$ \textit{%
with} $n\geq 2$ \textit{its finite-dimensional distributions are}%
\begin{eqnarray*}
&&\mathbb{P}_{\bar{\mu}}\left( \bar{Z}_{t_{1}}\in F_{1},...,\bar{Z}%
_{t_{n}}\in F_{n}\right) \\
&=&\sum_{\mathsf{m\in }\mathbb{N}^{d}}p_{\mathsf{m}}\int_{F_{1}}\mathsf{dx}%
_{1}...\int_{F_{n}}\mathsf{dx}_{n}\dprod\limits_{k=2}^{n}g\left( \mathsf{x}%
_{k},t_{k}-t_{k-1},\mathsf{x}_{k-1}\right) \\
&&\mathsf{\times }\left( \exp \left[ -t_{1}H\right] \varphi _{\mathsf{m,0}%
}\right) (\mathsf{x}_{1}\mathsf{)}\left( \exp \left[ -\left( T-t_{n}\right) H%
\right] \psi _{\mathsf{m,}T}\right) (\mathsf{x}_{n})
\end{eqnarray*}%
\textit{for all }$F_{1},...,F_{n}\in \mathcal{B}_{d}$\textit{\ and all }$%
0<t_{1}<...<t_{n}<T$\textit{. }

\textit{(b) We have}%
\begin{eqnarray*}
&&\mathbb{P}_{\bar{\mu}}\left( \bar{Z}_{t}\in F\right) \\
&=&\sum_{\mathsf{m\in }\mathbb{N}^{d}}p_{\mathsf{m}}\int_{F}\mathsf{dx}%
\left( \exp \left[ -tH\right] \varphi _{\mathsf{m,0}}\right) (\mathsf{x)}%
\left( \exp \left[ -\left( T-t\right) H\right] \psi _{\mathsf{m,}T}\right) (%
\mathsf{x})
\end{eqnarray*}%
\textit{for each }$F\in \mathcal{B}_{d}$\textit{\ and every} $t\in \left[ 0,T%
\right] $.

\textit{(c) We have}%
\begin{eqnarray}
&&\mathbb{E}_{_{\bar{\mu}}}\left( b(\bar{Z}_{t})\right)
\label{expectationter} \\
&=&\sum_{\mathsf{m\in }\mathbb{N}^{d}}p_{\mathsf{m}}\int_{\mathbb{R}^{d}}%
\mathsf{dx}b\left( \mathsf{x}\right) \left( \exp \left[ -tH\right] \varphi _{%
\mathsf{m,0}}\right) (\mathsf{x)}\left( \exp \left[ -\left( T-t\right) H%
\right] \psi _{\mathsf{m,}T}\right) (\mathsf{x})  \notag
\end{eqnarray}%
\textit{for each bounded Borel measurable function }$b:\mathbb{R}^{d}\mathbb{%
\mapsto C}$ \textit{and every }$t\in \left[ 0,T\right] $\textit{.}

\bigskip

\textbf{Proof.} The substitution of (\ref{density}) into (\ref{distribution}%
) gives%
\begin{eqnarray*}
&&\mathbb{P}_{\mu }\left( Z_{t_{1}}\in F_{1},...,Z_{t_{n}}\in F_{n}\right) \\
&=&\sum_{\mathsf{m\in }\mathbb{N}^{d}}p_{\mathsf{m}}\int_{F_{1}}\mathsf{dx}%
_{1}...\int_{F_{n}}\mathsf{dx}_{n}\dprod\limits_{k=2}^{n}g\left( \mathsf{x}%
_{k},t_{k}-t_{k-1},\mathsf{x}_{k-1}\right) \\
&&\times \int_{\mathbb{R}^{d}\times \mathbb{R}^{d}}\mathsf{d\mathsf{xdy}}%
\varphi _{\mathsf{m,0}}(\mathsf{x)}\psi _{\mathsf{m,}T}(\mathsf{y})g\left( 
\mathsf{x}_{1},t_{1},\mathsf{x}\right) g\left( \mathsf{y},T-t_{n},\mathsf{x}%
_{n}\right)
\end{eqnarray*}%
for all $F_{1},...,F_{n}\in \mathcal{B}_{d}$\ and all $0<t_{1}<...<t_{n}<T$,
where we have used the fact that $t_{0}=0$ and $\mathsf{x}_{0}\mathsf{=x}$.
This proves (a) since%
\begin{equation*}
\left( \exp \left[ -t_{1}H\right] \varphi _{\mathsf{m,0}}\right) (\mathsf{x}%
_{1}\mathsf{)=}\int_{\mathbb{R}^{d}}\mathsf{d\mathsf{x}}g\left( \mathsf{x}%
_{1},t_{1},\mathsf{x}\right) \varphi _{\mathsf{m,0}}(\mathsf{x)}
\end{equation*}%
and%
\begin{equation*}
\left( \exp \left[ -\left( T-t_{n}\right) H\right] \psi _{\mathsf{m,}%
T}\right) (\mathsf{x}_{n})=\int_{\mathbb{R}^{d}}\mathsf{\mathsf{dx}}g\left( 
\mathsf{x}_{n},T-t_{n},\mathsf{x}\right) \psi _{\mathsf{m,}T}(\mathsf{x}).
\end{equation*}%
The proof of (b) is similar by using (\ref{density}) in (\ref{probability1}%
). It is also plain that (c) follows from (b) and that $\bar{Z}_{\tau \in %
\left[ 0,T\right] }$ is non-stationary and non-Markovian for the same
reasons as those given in the proof of Theorem 2 of the preceding section. \
\ $\blacksquare $

\bigskip

When $\bar{\mu}$ is given by (\ref{positivedensity}) the associated process
remains stationary and it is useful to discuss its properties separately by
writing out the various quantities of interest in view of the applications
discussed in the next section:

\bigskip

\textbf{Corollary 1. }\textit{Assume that Hypothesis }(H)\textit{\ holds,
and let }$\bar{Z}_{\tau \in \left[ 0,T\right] }$ \textit{be the Bernstein
process in the sense of Proposition 1 with }$\bar{\mu}$\textit{\ given by (%
\ref{positivedensity}). Then the following statements are valid:}

\textit{(a) The process }$\bar{Z}_{\tau \in \left[ 0,T\right] }$ \textit{is
stationary, non-Markovian and for every }$n\in \mathbb{N}^{+}$ \textit{with} 
$n\geq 2$ \textit{its finite-dimensional distributions are}%
\begin{eqnarray}
&&\mathbb{P}_{\bar{\mu}}\left( \bar{Z}_{t_{1}}\in F_{1},...,\bar{Z}%
_{t_{n}}\in F_{n}\right)  \notag \\
&=&\mathcal{Z}^{-1}(T)\int_{F_{1}}\mathsf{dx}_{1}...\int_{F_{n}}\mathsf{dx}%
_{n}  \label{distributionquinto} \\
&&\times \dprod\limits_{k=2}^{n}g\left( \mathsf{x}_{k},t_{k}-t_{k-1},\mathsf{%
x}_{k-1}\right) \times g\left( \mathsf{x}_{1},T-(t_{n}-t_{1}),\mathsf{x}%
_{n}\right)  \notag
\end{eqnarray}%
\textit{for all }$F_{1},...,F_{n}\in \mathcal{B}_{d}$\textit{\ and all }$%
0<t_{1}<...<t_{n}<T$\textit{. }

\textit{(b) We have}%
\begin{equation}
\mathbb{P}_{\bar{\mu}}\left( \bar{Z}_{t}\in F\right) =\mathcal{Z}%
^{-1}(T)\int_{F}\mathsf{dx}g\left( \mathsf{x},T,\mathsf{x}\right)
\label{probability6}
\end{equation}%
\textit{for each }$F\in \mathcal{B}_{d}$\textit{\ and every} $t\in \left[ 0,T%
\right] $.

\textit{(c) We have}%
\begin{equation}
\mathbb{E}_{_{\bar{\mu}}}\left( b(\bar{Z}_{t})\right) =\mathcal{Z}%
^{-1}(T)\int_{\mathbb{R}^{d}}\mathsf{dx}b(\mathsf{x})g\left( \mathsf{x},T,%
\mathsf{x}\right)  \label{averageter}
\end{equation}%
\textit{for each bounded Borel measurable function }$b:\mathbb{R}^{d}\mathbb{%
\mapsto C}$ \textit{and every }$t\in \left[ 0,T\right] $\textit{.}

\bigskip

\textbf{Proof. }Relation (\ref{distributionquinto}) follows from the
substitution of (\ref{positivedensity}) into (\ref{distribution}) and from
the semigroup composition law for $g$, while (\ref{probability6}) is a
consequence of (\ref{positivedensity}) into (\ref{probability1}) and (\ref%
{averageter}) a consequence of (\ref{probability6}) since the density of the
law of the process is $\mathsf{x}\mapsto \mathcal{Z}^{-1}(T)g\left( \mathsf{x%
},T,\mathsf{x}\right) $. Now for any $\tau >0$ sufficiently small such that $%
0<t_{1}+\tau <...<t_{n}+\tau <T$ we have%
\begin{equation*}
\mathbb{P}_{\bar{\mu}}\left( \bar{Z}_{t_{1}+\tau }\in F_{1},...,\bar{Z}%
_{t_{n}+\tau }\in F_{n}\right) =\mathbb{P}_{\bar{\mu}}\left( \bar{Z}%
_{t_{1}}\in F_{1},...,\bar{Z}_{t_{n}}\in F_{n}\right)
\end{equation*}%
from (\ref{distributionquinto}) and therefore $\bar{Z}_{\tau \in \left[ 0,T%
\right] }$ is stationary, which entails the fact that both (\ref%
{probability6}) and (\ref{averageter}) are independent of $t$. Finally the
process is non-Markovian since (\ref{positivedensity}) is not of the form (%
\ref{markovianmeasure}). \ \ $\blacksquare $

\bigskip

As in the preceding section we can now define a linear transformation in $%
L^{2}\left( \mathbb{R}^{d}\right) $ which will play the r\^{o}le of a
density operator. Let us set%
\begin{equation}
\mathcal{R}(t)f:=\sum_{\mathsf{m\in }\mathbb{N}^{d}}p_{\mathsf{m}}\left(
f,u_{\mathsf{m}}(\mathsf{.},t)\right) _{2}v_{\mathsf{m}}(\mathsf{.},t)
\label{timedensityoperatorbis}
\end{equation}%
for each $f\in L^{2}\left( \mathbb{R}^{d}\right) $ and every $t\in \left[ 0,T%
\right] $, where $u_{\mathsf{m}}$ and $v_{\mathsf{m}}$ are given by (\ref%
{semigroup1}) and (\ref{semigroup2}), respectively. We begin with the
following:

\bigskip

\textbf{Lemma 3.} \textit{Assume that }$\left( \psi _{\mathsf{m,}T}\right) _{%
\mathsf{m\in }\mathbb{N}^{d}}$\textit{\ is an arbitrary bounded sequence in }%
$L^{2}\left( \mathbb{R}^{d}\right) $,\textit{\ and that }$\left( \varphi _{%
\mathsf{m},0}\right) _{\mathsf{m\in }\mathbb{N}^{d}}$\textit{\ is the
sequence associated with }$\left( \psi _{\mathsf{m,}T}\right) _{\mathsf{m\in 
}\mathbb{N}^{d}}$\textit{\ in the sense of Proposition 2. Then (\ref%
{timedensityoperatorbis}) defines a linear bounded operator in }$L^{2}\left( 
\mathbb{R}^{d}\right) $\textit{.}

\bigskip

\textbf{Proof. }Since the function $V$ is bounded from below according to
Hypothesis (H) we first have%
\begin{eqnarray}
\left\Vert u_{\mathsf{m}}(\mathsf{.},t)\right\Vert _{2} &\leq
&c_{T}\left\Vert \varphi _{\mathsf{m,0}}\right\Vert _{2},  \label{estimate8}
\\
\left\Vert v_{\mathsf{m}}(\mathsf{.},t)\right\Vert _{2} &\leq
&c_{T}\left\Vert \psi _{\mathsf{m,}T}\right\Vert _{2}  \label{estimate9}
\end{eqnarray}%
for some finite constant $c_{T}>0$ depending only on $T$ (and on the lower
bound in question). Moreover, by choosing $f=\varphi _{\mathsf{n,0}}$ in the
first inequality in (\ref{estimate2}) and by using the biorthogonality
relation (\ref{biorthogonality}), we see that for each $\theta \in \left[
0,1\right) $ there exists a finite constant $c_{\theta }>0$ such that $%
\left\Vert \varphi _{\mathsf{m,0}}\right\Vert _{2}\leq c_{\theta }$ for
every $\mathsf{m}\in \mathbb{N}^{d}$. Combining this with the boundedness of 
$\left( \psi _{\mathsf{m,}T}\right) _{\mathsf{m\in }\mathbb{N}^{d}}$ and
with (\ref{estimate8}), (\ref{estimate9}) we obtain%
\begin{eqnarray*}
\left\Vert u_{\mathsf{m}}(\mathsf{.},t)\right\Vert _{2} &\leq &c_{\ast
}<+\infty , \\
\left\Vert v_{\mathsf{m}}(\mathsf{.},t)\right\Vert _{2} &\leq &c_{\ast
}<+\infty
\end{eqnarray*}%
where $c_{\ast }$ depends only on $T$, the lower bound of $V$ and $\theta $.
Therefore we have%
\begin{equation*}
\sum_{\mathsf{m}\in \mathbb{N}^{d}}p_{\mathsf{m}}\left\vert \left( f,u_{%
\mathsf{m}}(\mathsf{.},t)\right) _{2}\right\vert \left\Vert v_{\mathsf{m}}(%
\mathsf{.},t)\right\Vert _{2}\leq c_{\ast }^{2}\left\Vert f\right\Vert
_{2}<+\infty
\end{equation*}%
for each $f\in L^{2}\left( \mathbb{R}^{d}\right) $ and every $t\in \left[ 0,T%
\right] $, so that (\ref{timedensityoperatorbis}) converges strongly in $%
L^{2}\left( \mathbb{R}^{d}\right) $ with the desired properties. \ \ $%
\blacksquare $

\bigskip

\textsc{Remark.} It is essential that the sequence $\left( \psi _{\mathsf{m,}%
T}\right) _{\mathsf{m\in }\mathbb{N}^{d}}$ be bounded for the above argument
to hold, but this does not follow from the first inequality in (\ref%
{estimate1}) as the boundedness of $\left( \varphi _{\mathsf{m},0}\right) _{%
\mathsf{m\in }\mathbb{N}^{d}}$ followed from the first inequality in (\ref%
{estimate2}). Indeed, the first inequality in (\ref{estimate1}) along with
the biorthogonality relation (\ref{biorthogonality}) only shows that there
exists a finite constant $c_{\theta }>0$ such that $\left\Vert \exp \left[
-TH\right] \psi _{\mathsf{m,}T}\right\Vert _{2}\leq c_{\theta }$ for every $%
\mathsf{m}\in \mathbb{N}^{d}$, but that does not entail the boundedness of $%
\left( \psi _{\mathsf{m,}T}\right) _{\mathsf{m\in }\mathbb{N}^{d}}$.

\bigskip

In fact we have much more than the conclusion of Lemma 3:

\bigskip

\textbf{Theorem 5. }\textit{The hypothesis is the same as in Lemma 3. Then
the following statements hold:}

\textit{(a) Expression (\ref{timedensityoperatorbis}) defines a linear
trace-class operator in }$L^{2}\left( \mathbb{R}^{d}\right) $\textit{\ with}%
\begin{eqnarray}
\func{Tr}\mathcal{R}(t) &=&\sum_{\mathsf{m\in }\mathbb{N}^{d}}p_{\mathsf{m}%
}=1,  \label{trace1} \\
\func{Tr}\mathcal{R}^{2}(t) &=&\sum_{\mathsf{m\in }\mathbb{N}^{d}}p_{\mathsf{%
m}}^{2}\leq 1  \label{trace2}
\end{eqnarray}%
\textit{for} \textit{every} $t\in \left[ 0,T\right] $. \textit{In particular
we have} $\func{Tr}\mathcal{R}^{2}(t)=1$\textit{\ if, and only if,} $p_{%
\mathsf{m}^{\ast }}=1$ \textit{for exactly one }$\mathsf{m}^{\ast }$, 
\textit{thus having }$p_{\mathsf{m}}=0$\textit{\ for every }$\mathsf{m\neq m}%
^{\ast }$\textit{.}

\textit{(b) The eigenvalue equations}%
\begin{eqnarray}
\mathcal{R}(t)v_{\mathsf{m}}(\mathsf{.},t) &=&p_{\mathsf{m}}v_{\mathsf{m}}(%
\mathsf{.},t),  \label{eigenequation1} \\
\mathcal{R}^{\ast }(t)u_{\mathsf{m}}(\mathsf{.},t) &=&p_{\mathsf{m}}u_{%
\mathsf{m}}(\mathsf{.},t)  \label{eigenequation2}
\end{eqnarray}%
\textit{hold for every }$\mathsf{m}\in \mathbb{N}^{d}$ \textit{and every} $%
t\in \left[ 0,T\right] $, \textit{where}%
\begin{equation}
\mathcal{R}^{\ast }(t)f=\sum_{\mathsf{m\in }\mathbb{N}^{d}}p_{\mathsf{m}%
}\left( f,v_{\mathsf{m}}(\mathsf{.},t)\right) _{2}u_{\mathsf{m}}(\mathsf{.}%
,t)  \label{adjointoperator}
\end{equation}%
\textit{is the adjoint of} $\mathcal{R}(t)$.

\textit{(c) Let us consider the linear bounded self-adjoint multiplication
operator on} $L^{2}\left( \mathbb{R}^{d}\right) $ \textit{given by} $Bf=bf$%
\textit{\ for every }$f\in L^{2}\left( \mathbb{R}^{d}\right) $,\textit{\
where }$b\in $\textit{\ }$L^{\infty }\left( \mathbb{R}^{d}\right) $\textit{\
is real-valued. If }$\bar{Z}_{\tau \in \left[ 0,T\right] }$ \textit{denotes
the Bernstein process of Theorem 4 then we have}%
\begin{equation}
\func{Tr}\left( \mathcal{R}\left( t\right) B\right) =\mathbb{E}_{_{\bar{\mu}%
}}\left( b(\bar{Z}_{t})\right)  \label{averagesequalitybis}
\end{equation}%
\textit{for every} $t\in \left[ 0,T\right] $, \textit{where the right-hand
side of (\ref{averagesequalitybis}) is given by (\ref{expectationter}).}

\bigskip

\textbf{Proof.} The proof of (\ref{trace1}) is quite similar to that of
Statement (a) in Theorem 3 and is thereby omitted, while that of (\ref%
{trace2}) follows from the biorthogonality of $u_{\mathsf{m}}(.,t)$ and $v_{%
\mathsf{m}}(.,t)$. Equations (\ref{eigenequation1}), (\ref{eigenequation2})
are an immediate consequence of (\ref{timedensityoperatorbis}), (\ref%
{adjointoperator}) and of the biorthogonality relation (\ref%
{biorthogonalitybis}), while the proof of (c) is identical to that of the
last statement of Theorem 3. \ \ $\blacksquare $

\bigskip

\textsc{Remark.} It follows directly from (\ref{biorthogonalitybis}) and (%
\ref{eigenequation1}) that%
\begin{equation*}
\sum_{\mathsf{m\in }\mathbb{N}^{d}}\left( \mathcal{R}\left( t\right) v_{%
\mathsf{m}}(\mathsf{.},t),u_{\mathsf{m}}(\mathsf{.},t)\right) _{2}=\sum_{%
\mathsf{m\in }\mathbb{N}^{d}}p_{\mathsf{m}}=1.
\end{equation*}%
Nevertheless, the fact that the preceding expression holds true is not
specific to the problem at hand, but is a general property of trace-class
operators whose trace may be computed by means of Lidskii's theorem using
biorthogonal systems generated by Riesz bases (see, e.g., Theorems 5 and 6
in Section 2, Chapter I, of \cite{gelfandvilenkin}). Finally, at the end of
this article we will dwell a bit more on the meaning of (\ref{trace2}).

\bigskip

If the initial-final conditions are given by (\ref{inifinacondbis}) we note
that (\ref{timedensityoperatorbis}) reduces to the self-adjoint, positive,
trace-class time-independent operator 
\begin{equation}
\mathcal{R}f=\sum_{\mathsf{m\in }\mathbb{N}^{d}}p_{\mathsf{m}}\left( f,%
\mathsf{f}_{\mathsf{m}}\right) _{2}\mathsf{f}_{\mathsf{m}}.
\label{densityoperator}
\end{equation}%
In this case we have the following:

\bigskip

\textbf{Corollary 2. }\textit{Let }$B:L^{2}\left( \mathbb{R}^{d}\right)
\mapsto L^{2}\left( \mathbb{R}^{d}\right) $ \textit{be the same operator as
in (c) of Theorem 5, and let }$\bar{Z}_{\tau \in \left[ 0,T\right] }$ 
\textit{be the Bernstein} \textit{process of Corollary 1. Then we have}%
\begin{equation}
\func{Tr}\left( \mathcal{R}B\right) \mathcal{=}\mathbb{E}_{\bar{\mu}%
}(b\left( \bar{Z}_{t}\right) )=\mathcal{Z}^{-1}(T)\sum_{\mathsf{m\in }%
\mathbb{N}^{d}}\exp \left[ -TE_{\mathsf{m}}\right] \int_{\mathbb{R}^{d}}%
\mathsf{dx}b\mathsf{(x)}\left\vert \mathsf{f}_{\mathsf{m}}\mathsf{(x)}%
\right\vert ^{2}  \label{averages}
\end{equation}%
\textit{for every }$t\in \left[ 0,T\right] $.

\bigskip

\textbf{Proof.} It is easily verified that%
\begin{equation*}
\func{Tr}\left( \mathcal{R}B\right) =\sum_{\mathsf{m\in }\mathbb{N}^{d}}p_{%
\mathsf{m}}\int_{\mathbb{R}^{d}}\mathsf{dx}b\mathsf{(x)}\left\vert \mathsf{f}%
_{\mathsf{m}}\mathsf{(x)}\right\vert ^{2}
\end{equation*}%
whenever $p_{\mathsf{m}}>0$ satisfies the normalization condition in (\ref%
{probabilities}), while if the probabilities associated with the spectrum
are given by (\ref{probagibbs}) we have%
\begin{eqnarray*}
&&\mathbb{E}_{\bar{\mu}}\left( b\left( \bar{Z}_{t}\right) \right) =\mathcal{Z%
}^{-1}(T)\int_{\mathbb{R}^{d}}\mathsf{dx}b\mathsf{(x)}g\left( \mathsf{x},T,%
\mathsf{x}\right) \\
&=&\mathcal{Z}^{-1}(T)\sum_{\mathsf{m\in }\mathbb{N}^{d}}\exp \left[ -TE_{%
\mathsf{m}}\right] \int_{\mathbb{R}^{d}}\mathsf{dx}b\mathsf{(x)}\left\vert 
\mathsf{f}_{\mathsf{m}}\mathsf{(x)}\right\vert ^{2}
\end{eqnarray*}%
for every $t\in \left[ 0,T\right] $ according to (\ref{expansion}) and (\ref%
{averageter}). \ \ $\blacksquare $

\bigskip

In the final section of this article we apply some of the above results to
the class of Bernstein processes generated by (\ref{cauchyforwardter}) and (%
\ref{cauchybackwardter}).

\section{On the periodic Ornstein-Uhlenbeck process and related processes}

We begin by recalling that the eigenvalue equation%
\begin{equation*}
\left( -\frac{1}{2}\Delta _{\mathsf{x}}+\frac{\lambda ^{2}}{2}\left\vert 
\mathsf{x}\right\vert ^{2}\right) \mathsf{h}_{\mathsf{m},\lambda }\mathsf{%
(x)=}E_{\mathsf{m},\lambda }\mathsf{h}_{\mathsf{m},\lambda }\mathsf{(x)}
\end{equation*}%
holds for every $\mathsf{m\in }$ $\mathbb{N}^{d}$, with%
\begin{equation}
E_{\mathsf{m},\lambda }:\mathsf{=}\left( \sum_{j=1}^{d}m_{j}+\frac{d}{2}%
\right) \lambda  \label{spectrum}
\end{equation}%
and%
\begin{equation*}
\mathsf{h}_{\mathsf{m},\lambda }\mathsf{(x):=}\dprod%
\limits_{j=1}^{d}h_{m_{j},\lambda }(x_{j}).
\end{equation*}%
In these expressions $m_{j}$ is the $j^{th}$ component of $\mathsf{m}$, $%
x_{j}$ the $j^{th}$ component of $\mathsf{x}$ and $h_{m,\lambda }$ denotes
the one-dimensional, suitably scaled Hermite function%
\begin{equation*}
h_{m,\lambda }(x):=\sqrt[4]{\lambda }h_{m}\left( \sqrt{\lambda }x\right)
\end{equation*}%
where%
\begin{equation*}
h_{m}\left( x\right) =(-1)^{m}\left( \pi ^{\frac{1}{2}}2^{m}m!\right) ^{-%
\frac{1}{2}}e^{\frac{x^{2}}{2}}\frac{d^{m}}{dx^{m}}e^{-x^{2}}\text{.}
\end{equation*}%
Furthermore we have%
\begin{equation}
\mathcal{Z}_{\lambda }(T):=\sum_{\mathsf{m}\in \mathbb{N}^{d}}\exp \left[
-TE_{\mathsf{m,}\lambda }\right] =\left( 2\left( \cosh (\lambda T)-1\right)
\right) ^{-\frac{d}{2}}  \label{convergencebis}
\end{equation}%
by summing the series explicitly, so that Mehler's kernel (\ref{mehler}) may
be expanded as%
\begin{equation*}
g_{\lambda }(\mathsf{x},t,\mathsf{y})=\sum_{\mathsf{m}\in \mathbb{N}%
^{d}}\exp \left[ -tE_{\mathsf{m,}\lambda }\right] \mathsf{h}_{\mathsf{m}%
,\lambda }(\mathsf{x})\mathsf{h}_{\mathsf{m},\lambda }(\mathsf{y})
\end{equation*}%
according to the considerations of Section 1, where the series is now
convergent for every $t\in \left( 0,T\right] $ uniformly in all $\mathsf{x},%
\mathsf{y\in }\mathbb{R}^{d}$. This last property is a consequence of the
Cram\'{e}r-Charlier inequality%
\begin{equation*}
\left\vert \mathsf{h}_{\mathsf{m},\lambda }(\mathsf{x})\mathsf{h}_{\mathsf{m}%
,\lambda }(\mathsf{y})\right\vert \mathsf{\leq }\left( \frac{\lambda }{\pi
^{2}}\right) ^{\frac{d}{4}}k^{2d}
\end{equation*}%
valid with $k\leq 1.086435$ uniformly in all $\mathsf{x},\mathsf{y}$ and $%
\mathsf{m}$ (see, e.g., Section 10.18 in \cite{erdmagobertri} for the
one-dimensional case). We first illustrate some of the consequences of
Theorem 1 by considering the initial-final data%
\begin{equation}
\left\{ 
\begin{array}{c}
\varphi _{\mathsf{m,0}}(\mathsf{x})=\mathcal{N}_{\mathsf{m},\lambda }\delta
\left( \mathsf{x}\right) , \\ 
\\ 
\psi _{\mathsf{m,}T}(\mathsf{x})=\mathcal{N}_{\mathsf{m},\lambda }\delta
\left( \mathsf{x}-\mathsf{b}_{\mathsf{m}}\right)%
\end{array}%
\right.  \label{inifinacondter}
\end{equation}%
where $\left( \mathsf{b}_{\mathsf{m}}\right) _{\mathsf{m}\in \mathbb{N}%
^{d}}\subset \mathbb{R}^{d}$ is an arbitrary sequence of points associated
with (\ref{spectrum}), and where%
\begin{equation*}
\mathcal{N}_{\mathsf{m},\lambda }:=\left( \frac{2\pi \sinh \left( \lambda
T\right) }{\lambda }\right) ^{\frac{d}{4}}\exp \left[ \frac{\lambda \coth
\left( \lambda T\right) \left\vert \mathsf{b}_{\mathsf{m}}\right\vert ^{2}}{4%
}\right] .
\end{equation*}%
A glance at (\ref{mehler}) shows that (\ref{inifinacondter}) is a particular
case of (\ref{inifinacond}) when $\mathsf{a}_{\mathsf{m}}=0$ for every $%
\mathsf{m}$. The corresponding solutions to (\ref{cauchyforwardter}) and (%
\ref{cauchybackwardter}) given by (\ref{forwardsolution bis}) and (\ref%
{backwardsolutionbis}) then read%
\begin{equation}
u_{\mathsf{m},\lambda }(\mathsf{x},t)=\mathcal{N}_{\mathsf{m},\lambda
}^{\ast }\sinh ^{-\frac{d}{2}}\left( \lambda t\right) \exp \left[ -\frac{%
\alpha _{\lambda }(t)\left\vert \mathsf{x}\right\vert ^{2}}{2}\right]
\label{cauchyforwardquarto}
\end{equation}%
and%
\begin{eqnarray}
v_{\mathsf{m},\lambda }(\mathsf{x},t) &=&\mathcal{N}_{\mathsf{m},\lambda
}^{\ast }\exp \left[ -\frac{\alpha _{\lambda }(T-t)\left\vert \mathsf{b}_{%
\mathsf{m}}\right\vert ^{2}}{2}\right] \sinh ^{-\frac{d}{2}}\left( \lambda
\left( T-t\right) \right)  \notag \\
&&\times \exp \left[ -\frac{1}{2}\left( \alpha _{\lambda }(T-t)\left\vert 
\mathsf{x}\right\vert ^{2}-\frac{2\lambda \left( \mathsf{b}_{\mathsf{m}},%
\mathsf{x}\right) _{\mathbb{R}^{d}}}{\sinh \left( \lambda \left( T-t\right)
\right) }\right) \right] ,  \label{cauchybackwardquarto}
\end{eqnarray}%
respectively, where we have defined%
\begin{equation}
\alpha _{\lambda }(t):=\lambda \coth \left( \lambda t\right)  \label{alpha}
\end{equation}%
for every $t\in \left( 0,T\right] $ and%
\begin{equation*}
\mathcal{N}_{\mathsf{m},\lambda }^{\ast }:=\left( \frac{\lambda \sinh \left(
\lambda T\right) }{2\pi }\right) ^{\frac{d}{4}}\exp \left[ \frac{\alpha
_{\lambda }(T)\left\vert \mathsf{b}_{\mathsf{m}}\right\vert ^{2}}{4}\right] .
\end{equation*}%
Then the following result holds:

\bigskip

\textbf{Corollary 3.}\textit{\ The Bernstein process }$Z_{\tau \in \left[ 0,T%
\right] }^{\mathsf{m},\lambda }$\textit{\ associated with (\ref%
{cauchyforwardter}), (\ref{cauchybackwardter}) and (\ref{inifinacondter}) in
the sense of Theorem 1 is a non-stationary Gaussian and Markovian process
such that the following properties are valid:}

\textit{(a) We have}%
\begin{equation}
\mathbb{P}_{\mu _{\mathsf{m}}}\left( Z_{t}^{\mathsf{m,}\lambda }\in F\right)
=(2\pi \sigma _{\lambda }(t))^{-\frac{d}{2}}\int_{F}d\mathsf{x}\exp \left[ -%
\frac{\left\vert \mathsf{x-b}_{\mathsf{m,}\lambda }(t)\right\vert ^{2}}{%
2\sigma _{\lambda }(t)}\right]  \label{probability11}
\end{equation}%
\textit{for each }$t\in \left( 0,T\right) $\textit{\ and every} $F\in 
\mathcal{B}_{d}$, \textit{where}%
\begin{equation}
\mathsf{b}_{\mathsf{m,}\lambda }(t)=\frac{\sinh (\lambda t)}{\sinh (\lambda
T)}\mathsf{b}_{\mathsf{m}}  \label{meanvector}
\end{equation}%
\textit{and}%
\begin{equation}
\sigma _{\lambda }(t)=\frac{\sinh \left( \lambda (T-t\right) )\sinh (\lambda
t)}{\lambda \sinh (\lambda T)}.  \label{variance}
\end{equation}%
\textit{Furthermore we have}%
\begin{equation}
\mathbb{P}_{\mu _{\mathsf{m}}}\left( Z_{0}^{\mathsf{m,}\lambda }=\mathsf{o}%
\right) =\mathbb{P}_{\mu _{\mathsf{m}}}\left( Z_{T}^{\mathsf{m,}\lambda }=%
\mathsf{b}_{\mathsf{m}}\right) =1  \label{probability12}
\end{equation}%
\textit{for every} $\mathsf{m}\in \mathbb{N}^{d}$.

\textit{(b) We have}%
\begin{equation}
\mathbb{E}_{\mu _{\mathsf{m}}}\left( (Z_{s}^{\mathsf{m,}\lambda ,i}-b_{%
\mathsf{m,}\lambda }^{i}(s))(Z_{t}^{\mathsf{m,}\lambda ,j}-b_{\mathsf{m}%
,\lambda }^{j}(t))\right) =\left\{ 
\begin{array}{c}
\frac{\sinh \left( \lambda (T-t\right) )\sinh (\lambda s)}{\lambda \sinh
(\lambda T)}\delta _{i,j}\text{ for }t\geq s, \\ 
\\ 
\frac{\sinh \left( \lambda (T-s\right) )\sinh (\lambda t)}{\lambda \sinh
(\lambda T)}\delta _{i,j}\text{ for }t\leq s%
\end{array}%
\right.  \label{covariance}
\end{equation}%
\textit{for all }$s,t\in \left[ 0,T\right] $\textit{\ and all }$i,j\in
\left\{ 1,...,d\right\} $\textit{, where }$b_{\mathsf{m,}\lambda }^{i}$%
\textit{\ denotes the }$i^{th}$\textit{\ component of }$\mathsf{b}_{\mathsf{%
m,}\lambda }$.

\textit{(c) We have}%
\begin{equation}
\mathbb{E}_{_{\mu _{\mathsf{m}}}}\left( b(Z_{t}^{\mathsf{m,}\lambda
})\right) =(2\pi \sigma _{\lambda }(t))^{-\frac{d}{2}}\int_{\mathbb{R}^{d}}%
\mathsf{dx}b(\mathsf{x})\exp \left[ -\frac{\left\vert \mathsf{x-b}_{\mathsf{%
m,}\lambda }(t)\right\vert ^{2}}{2\sigma _{\lambda }(t)}\right]
\label{averagequarto}
\end{equation}%
\textit{for each bounded Borel measurable function }$b:\mathbb{R}^{d}\mathbb{%
\mapsto C}$ \textit{and every }$t\in \left( 0,T\right) $.

\bigskip

\textbf{Proof.} We begin by proving (\ref{probability11}). Using (\ref%
{cauchyforwardquarto}) and (\ref{cauchybackwardquarto}) we first have%
\begin{eqnarray}
&&u_{\mathsf{m},\lambda }(\mathsf{x},t)v_{\mathsf{m},\lambda }(\mathsf{x},t)
\notag \\
&=&\left( \frac{\lambda \sinh (\lambda T)}{2\pi \sinh \left( \lambda
(T-t\right) )\sinh (\lambda t)}\right) ^{\frac{d}{2}}\exp \left[ \frac{%
\left( \alpha _{\lambda }(T)-\alpha _{\lambda }(T-t)\right) \left\vert 
\mathsf{b}_{\mathsf{m}}\right\vert ^{2}}{2}\right]  \notag \\
&&\times \exp \left[ -\frac{1}{2}\left( \left( \alpha _{\lambda }(t)+\alpha
_{\lambda }(T-t)\right) \left\vert \mathsf{x}\right\vert ^{2}-\frac{2\lambda
\left( \mathsf{b}_{\mathsf{m}},\mathsf{x}\right) _{\mathbb{R}^{d}}}{\sinh
\left( \lambda \left( T-t\right) \right) }\right) \right]
\label{probadensity}
\end{eqnarray}%
after regrouping terms, and furthermore%
\begin{eqnarray}
\alpha _{\lambda }(T)-\alpha _{\lambda }(T-t) &=&-\frac{\lambda \sinh
(\lambda t)}{\sinh \left( \lambda (T-t\right) )\sinh (\lambda T)}
\label{alphadiff} \\
\alpha _{\lambda }(t)+\alpha _{\lambda }(T-t) &=&\frac{\lambda \sinh
(\lambda T)}{\sinh \left( \lambda (T-t\right) )\sinh (\lambda t)}
\label{alphasum}
\end{eqnarray}%
from (\ref{alpha}). The substitution of (\ref{alphadiff}) and (\ref{alphasum}%
) into (\ref{probadensity}) then leads to%
\begin{eqnarray}
&&u_{\mathsf{m},\lambda }(\mathsf{x},t)v_{\mathsf{m,}\lambda }(\mathsf{x},t)
\notag \\
&=&\left( \frac{\lambda \sinh (\lambda T)}{2\pi \sinh \left( \lambda
(T-t\right) )\sinh (\lambda t)}\right) ^{\frac{d}{2}}\exp \left[ -\frac{%
\lambda \sinh (\lambda t)\left\vert \mathsf{b}_{\mathsf{m}}\right\vert ^{2}}{%
2\sinh \left( \lambda (T-t\right) )\sinh (\lambda T)}\right]  \notag \\
&&\times \exp \left[ -\frac{\lambda }{2}\left( \frac{\sinh (\lambda
T)\left\vert \mathsf{x}\right\vert ^{2}-2\sinh (\lambda t)\left( \mathsf{b}_{%
\mathsf{m}},\mathsf{x}\right) _{\mathbb{R}^{d}}}{\sinh \left( \lambda
(T-t\right) )\sinh (\lambda t)}\right) \right] .  \label{probadensityter}
\end{eqnarray}%
Now, for the numerator of the argument in the second exponential of the
preceding expression we have 
\begin{eqnarray}
&&\sinh (\lambda T)\left\vert \mathsf{x}\right\vert ^{2}-2\sinh (\lambda
t)\left( \mathsf{b}_{\mathsf{m}},\mathsf{x}\right) _{\mathbb{R}^{d}}  \notag
\\
&=&\sinh (\lambda T)\left\vert \mathsf{x-b}_{\mathsf{m,}\lambda
}(t)\right\vert ^{2}-\frac{\sinh ^{2}(\lambda t)\left\vert \mathsf{b}_{%
\mathsf{m}}\right\vert ^{2}}{\sinh (\lambda T)}  \label{identitybis}
\end{eqnarray}%
by virtue of (\ref{meanvector}). Therefore, taking (\ref{variance}) and (\ref%
{identitybis}) into account in (\ref{probadensityter}) we get%
\begin{eqnarray*}
&&u_{\mathsf{m},\lambda }(\mathsf{x},t)v_{\mathsf{m},\lambda }(\mathsf{x},t)
\\
&=&(2\pi \sigma _{\lambda }(t))^{-\frac{d}{2}}\exp \left[ -\frac{\left\vert 
\mathsf{x-b}_{\mathsf{m,}\lambda }(t)\right\vert ^{2}}{2\sigma _{\lambda }(t)%
}\right]
\end{eqnarray*}%
following the cancellation of two exponential factors, which proves
Statement (a) according to (\ref{probability5}). We also remark that (\ref%
{probability12}) is a particular case of (\ref{probability4}), and that (\ref%
{averagequarto}) holds according to (\ref{expectations}).

We now turn to the proof of (\ref{covariance}). According to (\ref%
{distributionoctavo}) we note that the density of the law of $(Z_{t_{1}}^{%
\mathsf{m,}\lambda },...,Z_{t_{n}}^{\mathsf{m,}\lambda })\in \mathbb{R}^{nd}$
is%
\begin{eqnarray*}
&&\dprod\limits_{k=2}^{n}g_{\lambda }\left( \mathsf{x}_{k},t_{k}-t_{k-1},%
\mathsf{x}_{k-1}\right) \times u_{\mathsf{m},\lambda }(\mathsf{x}%
_{1},t_{1})v_{\mathsf{m},\lambda }(\mathsf{x}_{n},t_{n}) \\
&=&\left( 2\pi \right) ^{-\frac{nd}{2}}\left( \frac{\lambda ^{n}\sinh
(\lambda T)}{\sinh \left( \lambda (T-t_{n}\right) )\sinh (\lambda t_{1})}%
\right) ^{\frac{d}{2}}\left( \dprod\limits_{k=2}^{n}\sinh \left( \lambda
(t_{k}-t_{k-1})\right) \right) ^{-\frac{d}{2}} \\
&&\times \exp \left[ \frac{1}{2}(\alpha _{\lambda }(T)-\alpha _{\lambda
}(T-t_{n}))\left\vert \mathsf{b}_{\mathsf{m}}\right\vert ^{2}\right] \\
&&\times \exp \left[ -\frac{\lambda }{2}\sum_{k=2}^{n}\frac{\cosh (\lambda
(t_{k}-t_{k-1}))\left( \left\vert \mathsf{x}_{k}\right\vert ^{2}+\left\vert 
\mathsf{x}_{k-1}\right\vert ^{2}\right) -2\left( \mathsf{x}_{k}\mathsf{,x}%
_{k-1}\right) _{\mathbb{R}^{d}}}{\sinh \left( \lambda (t_{k}-t_{k-1})\right) 
}\right] \\
&&\times \exp \left[ -\frac{1}{2}\left( \alpha _{\lambda }(t_{1})\left\vert 
\mathsf{x}_{1}\right\vert ^{2}+\alpha _{\lambda }(T-t_{n})\left\vert \mathsf{%
x}_{n}\right\vert ^{2}\right) \right] \times \exp \left[ \frac{\lambda
\left( \mathsf{b}_{\mathsf{m}},\mathsf{x}_{n}\right) _{\mathbb{R}^{d}}}{%
\sinh \left( \lambda (T-t_{n}\right) )}\right]
\end{eqnarray*}%
for every $n\geq 2$. Therefore, the tridiagonal matrix $C_{\lambda }^{-1}$
corresponding to the quadratic part when $d=1$ is identified as%
\begin{equation*}
C_{\lambda ,k,k}^{-1}=\left\{ 
\begin{array}{c}
\frac{\lambda \sinh (\lambda t_{2})}{\sinh \left( \lambda
(t_{2}-t_{1}\right) )\sinh (\lambda t_{1})}\text{ \ \ for }k=1, \\ 
\\ 
\frac{\lambda \sinh (\lambda (t_{k+1}-t_{k-1}))}{\sinh \left( \lambda
(t_{k+1}-t_{k}\right) )\sinh (\lambda (t_{k}-t_{k-1}))}\text{ \ \ for }%
k=2,...,n-1, \\ 
\\ 
\frac{\lambda \sinh (\lambda (T-t_{n-1}))}{\sinh \left( \lambda
(T-t_{n}\right) )\sinh (\lambda (t_{n}-t_{n-1}))}\text{ \ \ for }k=n%
\end{array}%
\right.
\end{equation*}%
(the second line not being there if $n=2)$, and 
\begin{equation*}
C_{\lambda ,k,k-1}^{-1}=C_{\lambda ,k-1,k}^{-1}=-\frac{\lambda }{\sinh
\left( \lambda \left\vert t_{k}-t_{k-1}\right\vert \right) }\text{ \ \ for }%
k=2,...,n.
\end{equation*}%
Consequently, inverting the matrix and using numerous identities among
hyperbolic functions we eventually get%
\begin{equation*}
C_{\lambda ,k,l}=\left\{ 
\begin{array}{c}
\frac{\sinh \left( \lambda (T-t_{k})\right) \sinh \left( \lambda
t_{l}\right) }{\lambda \sinh \left( \lambda T\right) }\text{ \ \ for }k\geq
l, \\ 
\\ 
\frac{\sinh \left( \lambda (T-t_{l})\right) \sinh \left( \lambda
t_{k}\right) }{\lambda \sinh \left( \lambda T\right) }\text{ \ \ for }k\leq
l,%
\end{array}%
\right.
\end{equation*}%
which leads to (\ref{covariance}) by standard arguments. \ \ $\blacksquare $

\bigskip

\textsc{Remarks.} (1) Corollary 3 thus describes a sequence of random curves
all pinned down at the origin when $t=0$ and at $\mathsf{b}_{\mathsf{m}}$
when $t=T$, with probability one. We also remark that the Gaussian law is
not centered unless $\mathsf{b}_{\mathsf{m}}=\mathsf{o}$, and that the
process is clearly non-stationary and Markovian since (\ref{probability11})
depends explicitly on time and (\ref{covariance}) factorizes as the product
of a function of $s$ times a function of $t$. Moreover, we note that the
curve $\sigma _{\lambda }:\left[ 0,T\right] \mapsto \mathbb{R}_{0}^{+}$
given by (\ref{variance}) is concave aside from satisfying $\sigma _{\lambda
}(0)=\sigma _{\lambda }(T)=0$, and that it takes on the maximal value at the
midpoint of the time interval, namely,%
\begin{equation*}
\sigma _{\lambda }\left( \frac{T}{2}\right) =\frac{\sinh ^{2}\left( \frac{%
\lambda T}{2}\right) }{\lambda \sinh (\lambda T)},
\end{equation*}%
thereby retaining the main features of a Brownian bridge. In fact, $Z_{\tau
\in \left[ 0,T\right] }^{\mathsf{m},\lambda }$ does reduce to a Brownian
bridge in the limit $\lambda \rightarrow 0_{+}$ since%
\begin{equation*}
\lim_{\lambda \rightarrow 0_{+}}\mathbb{E}_{\mu _{\mathsf{m}}}\left( (Z_{s}^{%
\mathsf{m},\lambda ,i}-b_{\mathsf{m},\lambda }^{i}(s))(Z_{t}^{\mathsf{m}%
,\lambda ,j}-b_{\mathsf{m},\lambda }^{j}(t))\right) =\left\{ 
\begin{array}{c}
\frac{(T-t)s}{T}\delta _{i,j}\text{ \ \ \ for }t\geq s, \\ 
\\ 
\frac{(T-s)t}{T}\delta _{i,j}\text{ \ \ for }t\leq s%
\end{array}%
\right.
\end{equation*}%
according to (\ref{covariance}).

(2) Relations (\ref{inifinacondter}) represent a very degenerate case of
Gaussian data. It would have been possible to replace (\ref{inifinacondter})
by choosing genuine Gaussian curves for both $\varphi _{\mathsf{m,0}}$ and $%
\psi _{\mathsf{m,}T}$, or by 
\begin{equation*}
\left\{ 
\begin{array}{c}
\varphi _{\mathsf{m,0}}(\mathsf{x})=\exp \left[ \frac{\lambda Td}{4}\right]
\delta \left( \mathsf{x}\right) , \\ 
\\ 
\psi _{\mathsf{m,}T}(\mathsf{x})=\exp \left[ \frac{\lambda Td}{4}\right]
\exp \left[ -\frac{\lambda \left\vert x\right\vert ^{2}}{2}\right]%
\end{array}%
\right.
\end{equation*}%
for every $\mathsf{m}$. In this case the corresponding Bernstein process
would have been a non-stationary, Markovian centered process $Z_{\tau \in %
\left[ 0,T\right] }^{\lambda }$ satisfying%
\begin{equation*}
\mathbb{P}_{\mu }\left( Z_{0}^{\lambda }=\mathsf{o}\right) =1,
\end{equation*}%
whose variance and covariance are given by%
\begin{equation*}
\left\{ 
\begin{array}{c}
\sigma _{\lambda }(t)=\frac{\sinh \left( \lambda t\right) \exp \left[
-\lambda t\right] }{\lambda }, \\ 
\\ 
\mathbb{E}\left( Z_{s}^{\lambda ,i}Z_{t}^{\lambda ,j}\right) =\frac{\exp %
\left[ -\lambda \left( s+t\right) \right] }{2\lambda }\left( \exp \left[
2\lambda \left( s\wedge t\right) \right] -1\right) \delta _{i,j}%
\end{array}%
\right.
\end{equation*}%
respectively, in other words a process identical in law with a $d$%
-dimensional Ornstein-Uhlenbeck process conditioned to start at the origin.
We omit the details of the computations that led to the above formulae,
which are quite similar to those carried out above.

\bigskip

Finally, we still have the following consequence of Theorem 3, where the
density operator is defined by%
\begin{equation*}
\mathcal{R}_{\lambda }\left( t\right) f:=\sum_{\mathsf{m}\in \mathbb{N}%
^{d}}p_{\mathsf{m}}\left( f,u_{\mathsf{m,}\lambda }(\mathsf{.},t)\right)
_{2}v_{\mathsf{m,}\lambda }(\mathsf{.},t)
\end{equation*}%
for each complex-valued $f\in L^{2}\left( \mathbb{R}^{d}\right) $ and every $%
t\in \left( 0,T\right) $, where $u_{\mathsf{m,}\lambda }(\mathsf{.},t)$ and $%
v_{\mathsf{m,}\lambda }(\mathsf{.},t)$ are given by (\ref%
{cauchyforwardquarto}) and (\ref{cauchybackwardquarto}), respectively:

\bigskip

\textbf{Corollary 4}. \textit{Let }$\bar{Z}_{\tau \in \left[ 0,T\right]
}^{\lambda }$\textit{\ be the Bernstein process in the sense of Proposition
1 corresponding to the joint probability density }%
\begin{equation*}
\bar{\mu}_{\lambda }(\mathsf{x,y})=g_{\lambda }\left( \mathsf{x},T,\mathsf{y}%
\right) \delta \left( \mathsf{x}\right) \sum_{\mathsf{m\in }\mathbb{N}%
^{d}}p_{\mathsf{m}}\mathcal{N}_{\mathsf{m},\lambda }^{2}\delta \left( 
\mathsf{y}-\mathsf{b}_{\mathsf{m}}\right)
\end{equation*}%
\textit{\ generated from (\ref{inifinacondter}), where }$\left( \mathsf{b}_{%
\mathsf{m}}\right) _{\mathsf{m}\in \mathbb{N}^{d}}\subset \mathbb{R}^{d}$%
\textit{\ is an arbitrary sequence such that}%
\begin{equation*}
\sup_{\mathsf{m}\in \mathbb{N}^{d}}\left\vert \mathsf{b}_{\mathsf{m}%
}\right\vert <+\infty .
\end{equation*}%
\textit{If }$B$\textit{\ is the multiplication operator of Theorem 3, then
we have}%
\begin{equation*}
\func{Tr}\left( \mathcal{R}_{\lambda }\left( t\right) B\right) =(2\pi \sigma
_{\lambda }(t))^{-\frac{d}{2}}\sum_{\mathsf{m\in }\mathbb{N}^{d}}p_{\mathsf{m%
}}\int_{\mathbb{R}^{d}}\mathsf{dx}b(\mathsf{x})\exp \left[ -\frac{\left\vert 
\mathsf{x-b}_{\mathsf{m,}\lambda }(t)\right\vert ^{2}}{2\sigma _{\lambda }(t)%
}\right]
\end{equation*}%
\textit{for each }$t\in \left( 0,T\right) $\textit{\ and every }$p_{\mathsf{m%
}}>0$\textit{\ satisfying the normalization condition in (\ref{probabilities}%
).}

\bigskip

The situation is quite different from that we just described if we consider
the hierarchy (\ref{cauchyforwardter}), (\ref{cauchybackwardter}) with the
initial-final data%
\begin{eqnarray}
\varphi _{\mathsf{m,0,\lambda }}(\mathsf{x}) &=&\mathsf{h}_{\mathsf{%
m,\lambda }}(\mathsf{x}),  \notag \\
\psi _{\mathsf{m,}T,\lambda }(\mathsf{x}) &=&\exp \left[ TE_{\mathsf{%
m,\lambda }}\right] \mathsf{h}_{\mathsf{m,\lambda }}(\mathsf{x})
\label{otherchoice}
\end{eqnarray}%
and with (\ref{probagibbs}) for the probabilities associated with each level
of the spectrum, thus having%
\begin{equation}
\mathcal{R}_{\lambda }f:=\mathcal{Z}_{\lambda }^{-1}(T)\sum_{\mathsf{m\in }%
\mathbb{N}^{d}}\exp \left[ -TE_{\mathsf{m,}\lambda }\right] \left( f,\mathsf{%
h}_{\mathsf{m},\lambda }\right) _{2}\mathsf{h}_{\mathsf{m},\lambda }
\label{densityoperatorbis}
\end{equation}%
for the density operator (\ref{densityoperator}). Then we have:

\bigskip

\textbf{Theorem 6. }\textit{For every }$\lambda >0$ \textit{the Bernstein
process }$\bar{Z}_{\tau \in \left[ 0,T\right] }^{\lambda }$\textit{\
associated with the infinite hierarchy (\ref{cauchyforwardter})-(\ref%
{cauchybackwardter}) in the sense of Corollary 1 is a stationary,
non-Markovian Gaussian process such that the following statements are valid:}

(\textit{a) We have}%
\begin{equation}
\mathbb{P}_{\bar{\mu}}\left( \bar{Z}_{t}^{\lambda }\in F\right) =(2\pi
\sigma _{\lambda })^{-\frac{d}{2}}\int_{F}\mathsf{dx}\exp \left[ -\frac{%
\left\vert \mathsf{x}\right\vert ^{2}}{2\sigma _{\lambda }}\right]
\label{gaussianprocesster}
\end{equation}%
\textit{for each }$t\in \left[ 0,T\right] $ \textit{and every} $F\in 
\mathcal{B}_{d}$, \textit{where}%
\begin{equation}
\sigma _{\lambda }=\frac{\sinh \left( \lambda T\right) }{2\lambda \left(
\cosh (\lambda T)-1\right) }.  \label{variancebis}
\end{equation}

\textit{(b) The components of }$\bar{Z}_{\tau \in \left[ 0,T\right]
}^{\lambda }$\textit{\ satisfy the relation}%
\begin{equation}
\mathbb{E}_{\bar{\mu}}\left( \bar{Z}_{s}^{\lambda ,i}\bar{Z}_{t}^{\lambda
,j}\right) =\frac{\cosh \left( \lambda \left( \left\vert t-s\right\vert -%
\frac{T}{2}\right) \right) }{2\lambda \sinh \left( \frac{\lambda T}{2}%
\right) }\delta _{i,j}  \label{covariancebis}
\end{equation}%
\textit{for all }$s,t\in \left[ 0,T\right] $\textit{\ and all }$i,j\in
\left\{ 1,...,d\right\} $\textit{.}

\textit{(c) For every linear bounded self-adjoint multiplication operator }$%
B $\textit{\ on} $L^{2}\left( \mathbb{R}^{d}\right) $ \textit{as defined in
Theorem 5 we have}%
\begin{equation*}
\func{Tr}\left( \mathcal{R}_{\lambda }B\right) \mathcal{=}\mathbb{E}_{\bar{%
\mu}}(b\left( \bar{Z}_{t}^{\lambda }\right) )=(2\pi \sigma _{\lambda })^{-%
\frac{d}{2}}\int_{\mathbb{R}^{d}}\mathsf{dx}b(\mathsf{x})\exp \left[ -\frac{%
\left\vert \mathsf{x}\right\vert ^{2}}{2\sigma _{\lambda }}\right]
\end{equation*}%
\textit{for every} $t\in \left[ 0,T\right] $.

\bigskip

\textbf{Proof.} The process $\bar{Z}_{\tau \in \left[ 0,T\right] }^{\lambda
} $ is Gaussian by virtue of (\ref{distributionquinto}) with Green's
function (\ref{mehler}). Furthermore we have 
\begin{equation*}
g_{\lambda }(\mathsf{x},T,\mathsf{x})=\left( \frac{\lambda }{2\pi \sinh
\left( \lambda T\right) }\right) ^{\frac{d}{2}}\exp \left[ -\frac{\lambda
\left( \cosh (\lambda T)-\mathsf{1}\right) \left\vert \mathsf{x}\right\vert
^{2}}{\sinh \left( \lambda T\right) }\right]
\end{equation*}%
so that (\ref{gaussianprocesster}) with (\ref{variancebis}) follows
immediately from (\ref{probability6}) and (\ref{convergencebis}). We now
turn to the proof of (\ref{covariancebis}) by determining the Gaussian
density of $\left( \bar{Z}_{t_{1}}^{\lambda },...,\bar{Z}_{t_{n}}^{\lambda
}\right) $ in $\mathbb{R}^{nd}$ for any $n\in \mathbb{N}^{+}$ by
substituting (\ref{mehler}) and (\ref{convergencebis}) into (\ref%
{distributionquinto}). We obtain%
\begin{eqnarray*}
&&\left( 2\left( \cosh (\lambda T)-1\right) \right) ^{\frac{d}{2}} \\
&&\times \dprod\limits_{k=2}^{n}g_{\lambda }\left( \mathsf{x}%
_{k},t_{k}-t_{k-1},\mathsf{x}_{k-1}\right) \times g_{\lambda }\left( \mathsf{%
x}_{1},T-(t_{n}-t_{1}),\mathsf{x}_{n}\right) \\
&=&\left( 2\pi \right) ^{-\frac{nd}{2}}\left( \frac{2\lambda ^{n}\left(
\cosh (\lambda T)-1\right) }{\sinh \left( \lambda (T-(t_{n}-t_{1})\right) )}%
\right) ^{\frac{d}{2}}\dprod\limits_{k=2}^{n}\left( \sinh (\lambda
(t_{k}-t_{k-1})\right) ^{-\frac{d}{2}} \\
&&\times \exp \left[ -\frac{\lambda }{2}\sum_{k=2}^{n}\frac{\cosh (\lambda
(t_{k}-t_{k-1}))\left( \left\vert \mathsf{x}_{k}\right\vert ^{2}+\left\vert 
\mathsf{x}_{k-1}\right\vert ^{2}\right) -2\left( \mathsf{x}_{k}\mathsf{,x}%
_{k-1}\right) _{\mathbb{R}^{d}}}{\sinh \left( \lambda (t_{k}-t_{k-1})\right) 
}\right] \\
&&\times \exp \left[ -\frac{\lambda }{2}\frac{\cosh (\lambda
(T-(t_{n}-t_{1})))\left( \left\vert \mathsf{x}_{1}\right\vert
^{2}+\left\vert \mathsf{x}_{n}\right\vert ^{2}\right) -2\left( \mathsf{x}_{1}%
\mathsf{,x}_{n}\right) _{\mathbb{R}^{d}}}{\sinh \left( \lambda
(T-(t_{n}-t_{1}))\right) }\right] .
\end{eqnarray*}%
For the sake of clarity we identify the inverse of the covariance matrix $%
C_{\lambda }$ by considering the case $n=2$ separately from the case $n\geq
3 $. For $n=2$ we obtain%
\begin{equation*}
C_{\lambda ,k,k}^{-1}=\frac{\lambda \sinh (\lambda T)}{\sinh \left( \lambda
(t_{2}-t_{1}\right) )\sinh \left( \lambda (T-(t_{2}-t_{1})\right) )}\text{ \
\ \ \ for }k=1,2
\end{equation*}%
\bigskip and%
\begin{equation*}
C_{\lambda ,2,1}^{-1}=C_{\lambda ,1,2}^{-1}=-\frac{\lambda }{\sinh \left(
\lambda \left\vert t_{2}-t_{1}\right\vert \right) }\text{ }-\frac{\lambda }{%
\sinh \left( \lambda (T-\left\vert t_{2}-t_{1}\right\vert \right) )},
\end{equation*}%
while for $n\geq 3$ we get%
\begin{equation*}
C_{\lambda ,k,k}^{-1}=\left\{ 
\begin{array}{c}
\frac{\lambda \sinh \left( \lambda (T-(t_{n}-t_{2})\right) )}{\sinh \left(
\lambda (t_{2}-t_{1}\right) )\sinh \left( \lambda (T-(t_{n}-t_{1})\right) )}%
\text{ \ \ \ \ for }k=1, \\ 
\\ 
\frac{\lambda \sinh \left( \lambda (t_{k+1}-t_{k-1})\right) }{\sinh \left(
\lambda (t_{k+1}-t_{k})\right) \sinh \left( \lambda (t_{k}-t_{k-1})\right) }%
\text{ \ \ \ \ for }k=2,...,n-1, \\ 
\\ 
\frac{\lambda \sinh \left( \lambda (T-(t_{n-1}-t_{1})\right) )}{\sinh \left(
\lambda (t_{n}-t_{n-1})\right) \sinh \left( \lambda (T-(t_{n}-t_{1})\right) )%
}\text{ \ \ \ \ for }k=n\text{ ,}%
\end{array}%
\right.
\end{equation*}%
\begin{equation*}
C_{\lambda ,k,k-1}^{-1}=C_{\lambda ,k-1,k}^{-1}=-\frac{\lambda }{\sinh
\left( \lambda \left\vert t_{k}-t_{k-1}\right\vert \right) }\text{ \ \ for }%
k=2,...,n,
\end{equation*}%
and 
\begin{equation*}
C_{\lambda ,n,1}^{-1}=C_{\lambda ,1,n}^{-1}=-\frac{\lambda }{\sinh \left(
\lambda (T-\left\vert t_{n}-t_{1}\right\vert \right) )},
\end{equation*}%
all the remaining matrix elements being zero. In both cases we then obtain
by inversion%
\begin{equation*}
C_{\lambda ,k,l}=\frac{\sinh \left( \lambda \left\vert
t_{k}-t_{l}\right\vert \right) -\sinh \left( \lambda (\left\vert
t_{k}-t_{l}\right\vert -T)\right) }{2\lambda \left( \cosh (\lambda
T)-1\right) }
\end{equation*}%
for all $k,l\in \left\{ 1,...,n\right\} $ or, equivalently,%
\begin{equation*}
C_{\lambda ,k,l}=\frac{\cosh \left( \lambda \left( \left\vert
t_{k}-t_{l}\right\vert -\frac{T}{2}\right) \right) }{2\lambda \sinh \left( 
\frac{\lambda T}{2}\right) },
\end{equation*}%
so that (\ref{covariancebis}) eventually follows. Finally, Statement (c)
follows from Corollary 2 by taking (\ref{gaussianprocesster}) into account.
Note that independently of the considerations of the preceding section a
glance at (\ref{gaussianprocesster}) and (\ref{covariancebis}) shows
directly that $\bar{Z}_{\tau \in \left[ 0,T\right] }^{\lambda }$ is
stationary, as well as non-Markovian since (\ref{covariancebis}) does not
factorize as the product of a function of $s$ times a function of $t$. \ \ $%
\blacksquare $

\bigskip

\textsc{Remarks.} (1) It turns out that the process of Theorem 6 identifies
in law with the $d$-dimensional periodic Ornstein-Uhlenbeck process. In
order to see this we define%
\begin{equation}
X_{t}:=\frac{e^{-\lambda t}}{1-e^{-\lambda T}}\int_{0}^{T}e^{-\lambda
(T-\tau )}\mathsf{dW}_{\tau }+\int_{0}^{t}e^{-\lambda (t-\tau )}\mathsf{dW}%
_{\tau },\text{ \ \ }t\in \left[ 0,T\right]  \label{periodicOU}
\end{equation}%
where $\mathsf{W}_{\tau \in \left[ 0,T\right] }$ is a given Wiener process
in $\mathbb{R}^{d}$, and where the integrals in (\ref{periodicOU}) are both
forward It\^{o} integrals. It is known from a particular case of Theorem 2.1
in \cite{kwakernaak}, or from a direct computation using the rules of It\^{o}
calculus (see also Section 5 in \cite{roellythieullen} for the case $d=1$),
that (\ref{periodicOU}) may be viewed as a way of writing the forward
Ornstein-Uhlenbeck integral equation with random periodic boundary conditions%
\begin{eqnarray}
X_{t} &=&e^{-\lambda t}X_{0}+\int_{0}^{t}e^{-\lambda (t-\tau )}\mathsf{dW}%
_{\tau },\text{ \ \ }t\in \left[ 0,T\right] ,  \notag \\
X_{0} &=&X_{T}  \label{itoOU}
\end{eqnarray}%
when $\mathbb{E}\left( X_{0}\right) =0$, whose covariance is precisely (\ref%
{covariancebis}). Therefore, our analysis shows that the periodic
Ornstein-Uhlenbeck process may be viewed as a very special example of a
stationary and non-Markovian Bernstein process. Incidentally, that process
happens to be quite relevant to the mathematical investigation of certain
quantum systems in equilibrium with a thermal bath when the inverse
temperature is interpreted as the period. This is indeed a consequence of
the fact that it also identifies in law with the Gaussian process of mean
zero used in Theorem 2.1 of \cite{hoeghkrohn} when the positive matrix
therein is chosen as $A=\lambda \mathbb{I}_{d}$ with $\mathbb{I}_{d}$ the
identity in $\mathbb{R}^{d}$. This, in turn, follows immediately from (\ref%
{covariancebis}) which may also be written as 
\begin{equation*}
\mathbb{E}_{\bar{\mu}}\left( \bar{Z}_{s}^{\lambda ,i}\bar{Z}_{t}^{\lambda
,j}\right) =\frac{\exp \left[ -\lambda \left\vert t-s\right\vert \right]
+\exp \left[ -\lambda \left( T-\left\vert t-s\right\vert \right) \right] }{%
2\lambda \left( 1-\exp \left[ -\lambda T\right] \right) }\delta _{i,j}
\end{equation*}%
by virtue of the identity 
\begin{equation*}
\frac{\cosh \left( \lambda \left( t-\frac{T}{2}\right) \right) }{\sinh
\left( \frac{\lambda T}{2}\right) }=\frac{\exp \left[ -\lambda t\right]
+\exp \left[ -\lambda \left( T-t\right) \right] }{1-\exp \left[ -\lambda T%
\right] }
\end{equation*}%
valid for every $t\in \left[ 0,T\right] $. Finally, we observe that the
definition of a periodic process "indexed by the circle" that satisfies the
"two-sided Markov property on the circle" given in Section 4 of \cite%
{kleinlandau} is a very special case of our definition of a Bernstein
process given at the beginning of this paper. Indeed, a standard argument
shows that Relation (\ref{condiexpectations}) is equivalent to the statement
that $\mathcal{F}_{s}^{+}\vee $ $\mathcal{F}_{t}^{-}$ is conditionally
independent of the $\sigma $-algebra%
\begin{equation*}
\mathcal{F}_{\left[ s,t\right] }:=\sigma \left\{ Z_{\tau }^{-1}\left(
F\right) :\tau \in \left[ s,t\right] ,\text{ }F\in \mathcal{B}_{d}\right\}
\end{equation*}%
when%
\begin{equation*}
\mathcal{F}_{\left\{ s,t\right\} }=\sigma \left\{ Z_{s}^{-1}\left( F\right) ,%
\text{ }Z_{t}^{-1}\left( F\right) :F\in \mathcal{B}_{d}\right\}
\end{equation*}%
is given. In this respect we also refer the reader to \cite{carmatheo} and 
\cite{jamisonbis} for the stationary Gaussian case when $d=1$. More
generally, we remark that Problem (\ref{itoOU}) falls into the realm of a
much more general class of periodic linear stochastic differential equations
which were investigated by several authors, including \cite{kwakernaak}
where some of the multidimensional time-periodic processes considered there
were useful regarding the resolution of filtering, smoothing and prediction
problems.

(2) We complete this article with an observation concerning the
interpretation of (\ref{trace2}). Following the analogy with Quantum
Statistical Mechanics, we may say that the operator $\mathcal{R}(t)$
represents a so-called pure state when $\func{Tr}\mathcal{R}^{2}(t)=1$ and a
mixed state when $\func{Tr}\mathcal{R}^{2}(t)<1$ (see, e.g., \cite%
{vonneumann} for explanations regarding this terminology). In view of the
first part of Theorem 5, it is therefore legitimate to say that the
non-Markovian Bernstein processes that we constructed from the method of
Section 4 correspond to mixed states in the above sense. Similar
considerations hold for operator (\ref{densityoperatorbis}), which satisfies
the inequalities%
\begin{equation*}
0\leq \mathcal{R}_{\lambda }^{2}\leq \mathcal{R}_{\lambda }\leq \mathbb{I}
\end{equation*}%
in the sense of quadratic forms since $\mathcal{R}_{\lambda }$ is
self-adjoint, where $\mathbb{I}$ stands for the identity in $L^{2}\left( 
\mathbb{R}^{d}\right) $. In this case we always have $\func{Tr}\mathcal{R}%
_{\lambda }^{2}(t)<1$ by virtue of (\ref{probagibbs}), and the only process
that would correspond to a pure state in this context is the Markovian
process generated by the measure%
\begin{equation*}
\mu _{0,\lambda }(G)=\int_{G}\mathsf{dxdy}\varphi _{0\mathsf{,0,\lambda }}(%
\mathsf{x})g_{\lambda }(\mathsf{x},T,\mathsf{y})\psi _{0\mathsf{,}T\mathsf{%
,\lambda }}(\mathsf{y})
\end{equation*}%
which is of the form (\ref{markovianmeasure}), where $g_{\lambda }$ is
Mehler's kernel (\ref{mehler}) and $\varphi _{0\mathsf{,0,\lambda }}$, $\psi
_{0\mathsf{,}T\mathsf{,\lambda }}$ are given by (\ref{otherchoice}) with $%
\mathsf{m}=0$. There are, however, many other interesting Markovian
Bernstein processes associated with (\ref{cauchyforwardter})-(\ref%
{cauchybackwardter}) (see, e.g., \cite{zambrinibis}).

\bigskip

\textbf{Acknowledgements.} The research of both authors was supported by the
FCT of the Portuguese government under grant PTDC/MAT-STA/0975/2014. The
first author is also indebted to Marco Dozzi for stimulating discussions
regarding certain probabilistic aspects of this article.


\begin{thebibliography}{99}
\bibitem{albeyaza} \textsc{Albeverio, S., Yasue, K., Zambrini, J-C.,} 
\textit{Euclidean quantum mechanics: analytical approach,} Annales de
l'Institut Henri Poincar\'{e}, Physique Th\'{e}orique, \textbf{49} (1989)
259-308.

\bibitem{apostol} \textsc{Apostol, T. M., }\textit{Mathematical Analysis,}
Addison-Wesley Series in Mathematics, Addison-Wesley Publishing Company,
Inc., Reading (1974).

\bibitem{aronsonbis} \textsc{Aronson D. G.,} \textit{Bounds for the
fundamental solution of a parabolic equation,} Bulletin of the American
Mathematical Society \textbf{73} (1967) 890-896.

\bibitem{aronson} \textsc{Aronson D. G.,}\textit{\ Non-negative solutions of
linear parabolic equations,} Annali della Scuola Normale Superiore di Pisa 
\textbf{22} (1968) 607-694.

\bibitem{bernstein} \textsc{Bernstein, S., }\textit{Sur les liaisons entre
les grandeurs al\'{e}atoires}, in: Verhandlungen des Internationalen
Mathematikerkongress \textbf{1} (1932) 288-309.

\bibitem{beurling} \textsc{Beurling, A., }\textit{An automorphism of product
measures,} Annals of Mathematics \textbf{72} (1960) 189-200.

\bibitem{carlen} \textsc{Carlen, E. A.,} \textit{Conservative diffusions},
Communications in Mathematical Physics \textbf{94} (1984) 293-315.

\bibitem{carmatheo} \textsc{Carmichael, J. P., Masse, J. C., Theodorescu, R.,%
} \textit{Processus Gaussiens stationnaires r\'{e}ciproques sur un
intervalle,} Comptes Rendus de l'Acad\'{e}mie des Sciences, S\'{e}rie I, 
\textbf{295} (1982) 291-294.

\bibitem{cruzeirozambrini} \textsc{Cruzeiro, A. B., Zambrini, J.-C.,}\textit{%
\ Malliavin calculus and Euclidean quantum mechanics, I. Functional calculus,%
} Journal of Functional Analysis \textbf{96} (1991) 62-95.

\bibitem{erdmagobertri} \textsc{Erd\'{e}lyi, A., Magnus, W., Oberhettinger,
F., Tricomi, F. G.,} \textit{Higher Transcendental Functions,} \textit{II,}
McGraw-Hill, Inc., New York (1953).

\bibitem{galichon} \textsc{Galichon, A., }\textit{Optimal Transport Methods
in Economics}, Princeton University Press, Princeton (2016).

\bibitem{gelfandvilenkin} \textsc{Gel'fand, I. M., Vilenkin, N. Y.,} \textit{%
Generalized Functions IV: Applications of Harmonic Analysis,} Academic
Press, New York (1964).

\bibitem{gohbergkrein} \textsc{Gohberg, I. C., Krein, M. G.,} \textit{%
Introduction to the Theory of Linear Nonselfadjoint Operators in Hilbert
Space,} Translations of Mathematical Monographs \textbf{18}, American
Mathematical Society, Providence (1969).

\bibitem{hoeghkrohn} \textsc{H}$\oslash $\textsc{egh-Krohn, R,}\textit{\
Relativistic quantum statistical mechanics in two-dimensional space-time, }%
Communications in Mathematical Physics\textbf{\ 38} (1974) 195-224.

\bibitem{jamisonbis} \textsc{Jamison, B.,} \textit{Reciprocal processes: the
stationary Gaussian case,} The Annals of Mathematical Statistics \textbf{41}
(1970) 1624-1630.

\bibitem{jamison} \textsc{Jamison, B.,} \textit{Reciprocal processes,}
Zeitschrift f\"{u}r Wahrscheinlichkeitstheorie und Verwandte Gebiete \textbf{%
30} (1974) 65-86.

\bibitem{kato} \textsc{Kato, T.,} \textit{Perturbation Theory for Linear
Operators,} Grundlehren der mathematischen Wissenschaften \textbf{132},
Springer Verlag, New York (1984).

\bibitem{kleinlandau} \textsc{Klein, A., Landau, L. J.,} \textit{Periodic
Gaussian Osterwalder-Schrader positive} \textit{processes and the two-sided
Markov property on the circle,} Pacific Journal of Mathematics \textbf{94}
(1981) 341-367.

\bibitem{kwakernaak} \textsc{Kwakernaak, H.,} \textit{Periodic linear
differential stochastic processes,} SIAM Journal on Control \textbf{13}
(1975) 400-413.

\bibitem{lasalle} \textsc{Lassalle, R.,} \textit{Causal transport plans and
their Monge-Kantorovich problems,} Stochastic Analysis and Applications, 
\textbf{36} (2018), in press.

\bibitem{leonard} \textsc{L\'{e}onard, C.,} \textit{A survey of the Schr\"{o}%
dinger problem and some of its connections with optimal transport,} Discrete
and Continuous Dynamical Systems, Series A, \textbf{34} (2014) 1533-1574.

\bibitem{leoroezamb} \textsc{L\'{e}onard, C., Roelly, S., Zambrini, J.-C.}, 
\textit{Reciprocal processes. A measure-theoretical point of view,}
Probability Surveys \textbf{11} (2014) 237-269.

\bibitem{paleywiener} \textsc{Paley, R. E. A. C, Wiener, N.,}\textit{\
Fourier Transforms in the Complex Domain, }American Mathematical Society
Colloquium Publications XIX, American Mathematical Society, New York (1934).

\bibitem{pedersen} \textsc{Pedersen, J., }\textit{Periodic
Ornstein-Uhlenbeck processes driven by L\'{e}vy processes, }Journal of
Applied Probability \textbf{39}\textit{\ }(2002) 748-763.

\bibitem{reedsimon} \textsc{Reed, M., Simon, B.,} \textit{Methods of Modern
Mathematical Physics IV: Analysis of Operators,} Academic Press, New York
(1978).

\bibitem{riesznagy} \textsc{Riesz, F., Nagy, B. SZ.,} \textit{Functional
Analysis,} Dover Books on Mathematics, Dover (1990).

\bibitem{roellythieullen} \textsc{Roelly, S., Thieullen, M.,} \textit{A
characterization of reciprocal processes via an integration by parts formula
on the path space,} Probability Theory and Related Fields \textbf{123}
(2002) 97-120.

\bibitem{schroedinger} \textsc{Schr\"{o}dinger, E.,} \textit{Sur la th\'{e}%
orie relativiste de l'\'{e}lectron et l'interpr\'{e}tation de la m\'{e}%
canique quantique,} Annales de l'Institut Henri Poincar\'{e} \textbf{2}
(1932) 269-310.

\bibitem{villani} \textsc{Villani, C.,}\textit{\ Optimal Transport: Old and
New,} Grundlehren der Mathematischen Wissenschaften \textbf{338}, Springer,
New York (2009).

\bibitem{vonneumann} \textsc{von Neumann, J.,} \textit{Mathematical
Foundations of Quantum Mechanics,} Princeton Landmarks in Mathematics
Series, Princeton University Press (1996).

\bibitem{vuillermot} \textsc{Vuillermot, P.-A.,} \textit{On the time
evolution of Bernstein processes associated with a class of parabolic
equations, }Discrete and Continuous Dynamical Systems, Series B, \textbf{23}
(2018), in press.

\bibitem{vuizambrini} \textsc{Vuillermot, P.-A., Zambrini, J.-C.,} \textit{%
Bernstein diffusions for a class of linear parabolic partial differential
equations,} Journal of Theoretical Probability \textbf{27} (2014) 449-492.

\bibitem{zambrini} \textsc{Zambrini, J.-C.,} \textit{Variational processes
and stochastic versions of mechanics,} Journal of Mathematical Physics 
\textbf{27} (1986) 2307-2330.

\bibitem{zambrinibis} \textsc{Zambrini, J.-C.,} \textit{The research program
of Stochastic Deformation (with a view toward Geometric Mechanics), }in:%
\textit{\ }Stochastic Analysis: a Series of Lectures, Birkh\"{a}user
Progress in Probability book series \textbf{68, }Eds. R. Dalang, M. Dozzi,
F. Flandoli, F. Russo (2015), pp. 359-393.
\end{thebibliography}
\end{document}